\documentclass[11pt]{article}
\usepackage{amsmath,amssymb,amstext,dsfont,fancyvrb,float,fontenc,graphicx,subfigure, theorem}

\title{\Large\bf An optimality result about sample path properties of Operator Scaling Gaussian Random Fields}
\author{\sc M. Clausel and B. Vedel}
\date{}

\cleardoublepage 
\usepackage[strict]{changepage}
\def\abstractname{Abstract -}   
\def\abstract{\begin{adjustwidth}{1cm}{1cm} \par    \footnotesize \noindent {\bf \abstractname}
\def\endabstract{ \end{adjustwidth} \smallskip }}

{\theorembodyfont{\itshape}\newtheorem{theorem}{Theorem}[section]}
{\theorembodyfont{\itshape}\newtheorem{proposition}{Proposition}[section]}
{\theorembodyfont{\itshape}\newtheorem{definition}{Definition}[section]}
{\theorembodyfont{\itshape}\newtheorem{lemma}{Lemma}[section]}
{\theorembodyfont{\itshape}\newtheorem{corollary}{Corollary}[section]}
{\theorembodyfont{\rm}\newtheorem{notation}{Notation}[section]}
{\theorembodyfont{\rm}\newtheorem{remark}{Remark}[section]}
{\theorembodyfont{\rm}}
{\theorembodyfont{\rm }}

\def\rmd{\mathrm{d}}
\def\rme{\mathrm{e}}
\def\rmi{\mathrm{i}}

\begin{document}
\maketitle
\vskip 1.5em
 \begin{abstract}
We study the sample paths properties of Operator
scaling Gaussian random fields. Such fields are anisotropic
generalizations of anisotropic self-similar random fields as anisotropic Fractional Brownian Motion. Some characteristic properties of the anisotropy are revealed by the regularity of the sample paths. The sharpest way of measuring smoothness is related to these
anisotropies and thus to the geometry of these fields.
 \end{abstract}

\begin{keywords}
Operator scaling Gaussian random field, anisotropy, sample paths
properties, anisotropic Besov spaces
\end{keywords}

\begin{MSC}
60G15 60G18 60G60 60G17 42C40 46E35
\end{MSC}


\section{Introduction and motivations}
Random fields are now used for modeling in a wide
range of scientific areas including physics, engineering, hydrology, biology, economics and finance (see~\cite{30} and its bibliography). An important requirement is that the data thus modelled present strong anisotropies which therefore have to be present in the model. Many anisotropic random fields have therefore been proposed as natural models in various areas such as image processing, hydrology,
geostatistics and spatial statistics (see, for example, Davies and Hall~\cite{13}, Bonami and Estrade~\cite{7}, Benson et al.~\cite{3}). Let us also quote the example of Levy random fields, deeply studied by Durand and Jaffard (see~\cite{15}), which is the only known model of anisotropic multifractal random field. In many cases, Gaussian models have turned to be relevant when investigating anisotropic problems. For example the stochastic model of surface waves is usually assumed to be Gaussian and is surprisingly accurate (see~\cite{20}). More generally anisotropic Gaussian random fields are involved in many others concrete situations and then arise naturally in stochastic partial differential equations (see, e.g., Dalang~\cite{12}, Mueller and Tribe~\cite{23}, \^{O}ksendal and Zhang~\cite{26}, Nualart~\cite{25}).

In many situations, the data present invariant features across the scales (see for example~\cite{1}). These two requirements (anisotropy and self--similarity) may seem contradictory, since the classical notion of self--similarity defined for a random field $\{X(x)\}_{x\in\mathbb{R}^{d}}$ on $\mathbb{R}^d$ by
\begin{equation}\label{EqSS}
\{X(ax)\}_{x \in \mathbb{R}^d} \overset{\mathcal{L}}{=}  \{a^{H_0} X(x)\}_{x\in \mathbb{R}^d}\;,
\end{equation}
for some $H_0 \in \mathbb{R}$ (called the Hurst index) is by construction isotropic and has then to be changed in order to fit anisotropic situations. To this end, several extensions of self--similarity property in an anisotropic setting have been proposed. In~\cite{18}, Hudson and Mason defined operator self-similar processes $\{X(t)\}_{t\in\mathbb{R}}$ with values in $\mathbb{R}^d$. In~\cite{19}, Kamont introduced Fractional Brownian Sheets which satisfies different scaling properties according to the coordinate axes. More recently, in \cite{6} Bierm\'e, Meerschaert and Scheffler introduced the notion of Operator Scaling Random Fields (OSRF). These fields satisfy the following anisotropic scaling relation :
\begin{equation}\label{selfsimilar}
\{X(a^{E_0}x)\}_{x \in \mathbb{R}^d} \overset{\mathcal{L}}{=}  \{a^{H_0} X(x)\}_{x\in \mathbb{R}^d}\;,
\end{equation}
for some matrix $E_0$ (called an exponent or an anisotropy of the field) whose
eigenvalues have a positive real part and some $H_0>0$ (called an Hurst index of the field).
The usual notion of self-similarity is extended replacing {\bf
usual scaling}, (corresponding to the case $E_0=Id$)  by {\bf a
linear scaling} involving the matrix $E_0$ (see figure~\ref{Fig1}
below). It allows to define new classes of random fields with new
geometry and structure.
\begin{figure}[H]\label{Fig1}
\begin{minipage}[c]{.32\linewidth}
    \centering
 \includegraphics[angle=0,width=.95\textwidth]{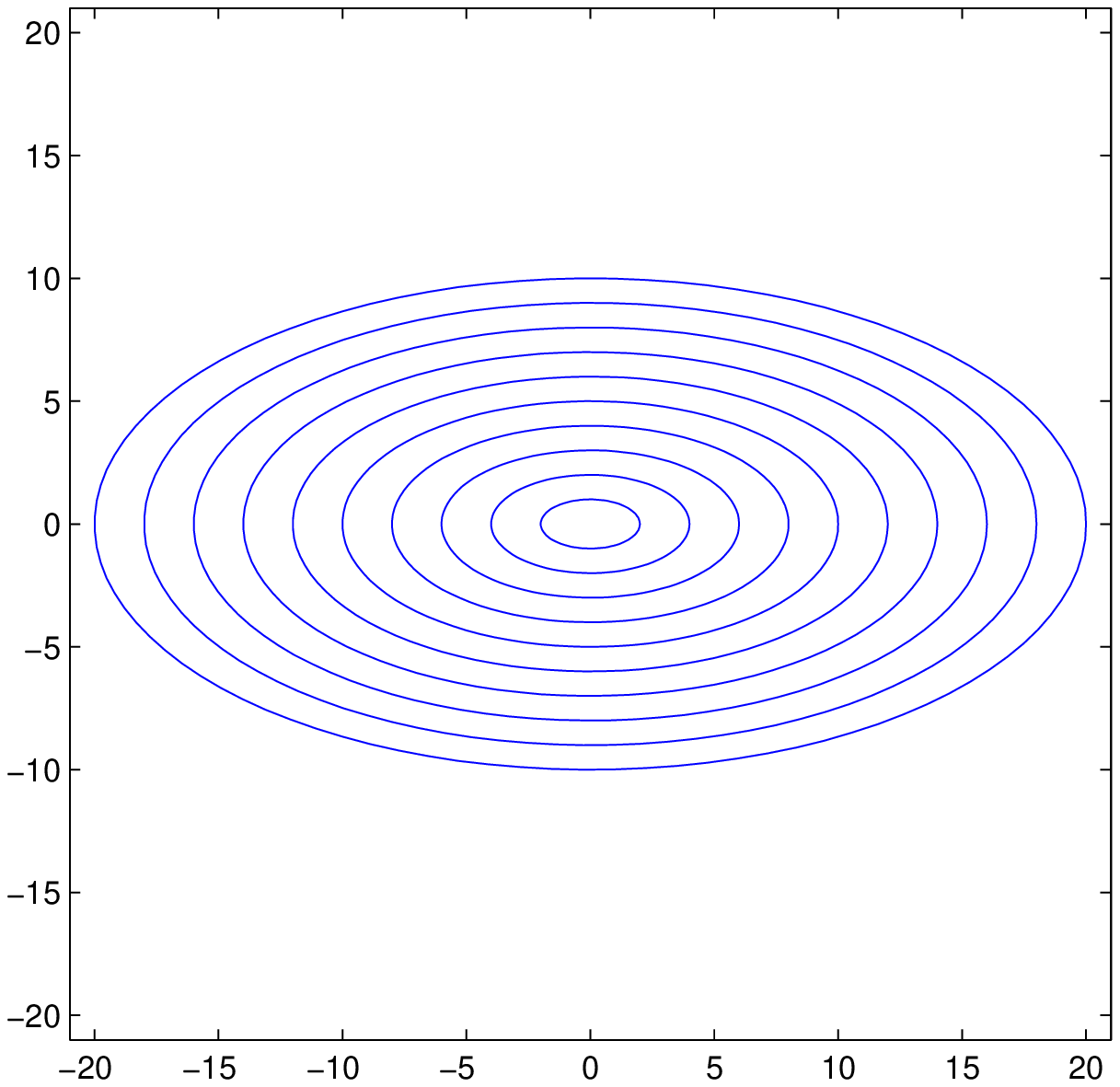}
\begin{center}$E=\begin{pmatrix}1&0\\0&1\end{pmatrix}$, $\lambda\in\{1,\cdots,10\}$\end{center}
\end{minipage}\hfill
  \begin{minipage}[c]{.32\linewidth}
    \hspace{.025\textwidth}
      \includegraphics[angle=0,width=.95\textwidth]{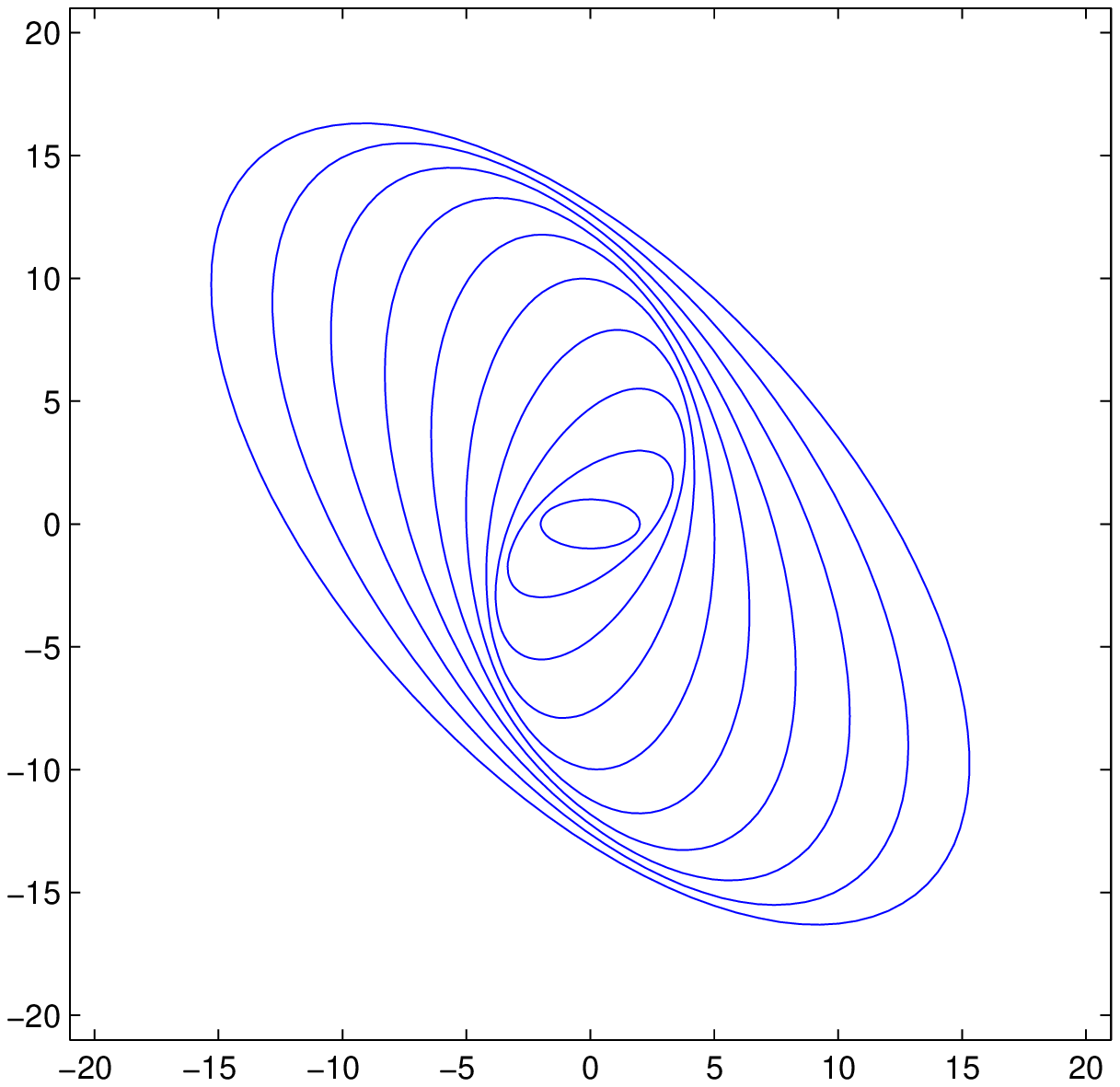}
      \begin{center}$E=\begin{pmatrix}1&-1\\1&1\end{pmatrix}$, $\lambda\in\{1,\cdots,10\}$\end{center}
\end{minipage}
 \begin{minipage}[c]{.32\linewidth}
    \hspace{.025\textwidth}
       \includegraphics[angle=0,width=.95\textwidth]{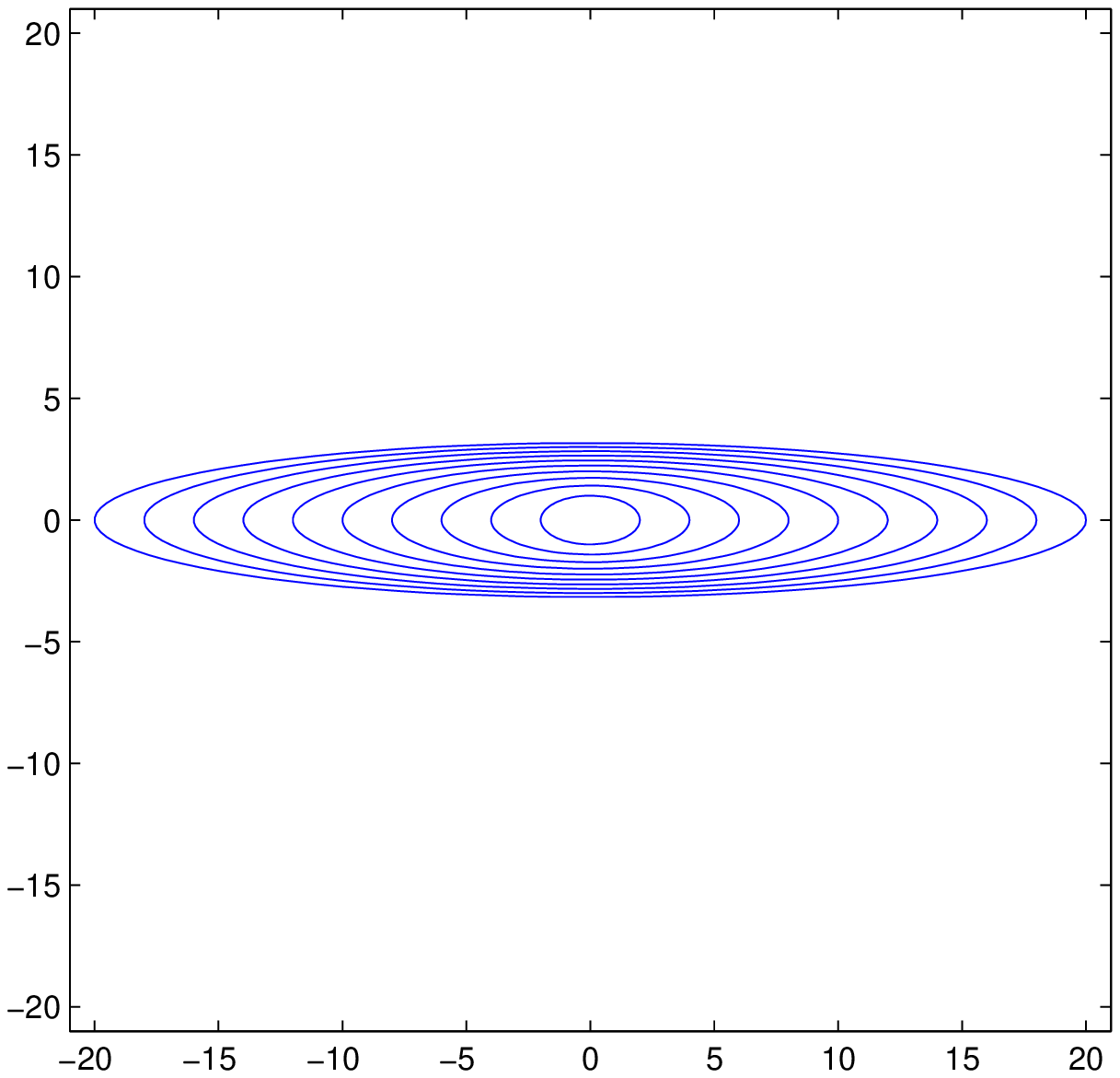}
       \begin{center}$E=\begin{pmatrix}1&0\\0&1/2\end{pmatrix}$, $\lambda\in\{1,\cdots,10\}$\end{center}
\end{minipage}
\begin{center}
\begin{caption} \;Action of a linear scaling $x\mapsto \lambda^{E} x$ on the smallest ellipsis.\end{caption}
\end{center}
\end{figure}
This new class of random fields have been introduced in order to model various phenomena such as fracture surfaces (see~\cite{27}) or sedimentary aquifers (see~\cite{3}). In~\cite{6}, the authors construct a large class of Operator Scaling Stable Random
Fields with stationary increments presenting both a moving average and an harmonizable representation of these fields.

In order to use such models in practice, the first problem is to recover the parameters $H_0$ and $E_0$ from the inspection of one sample paths. Even if we consider the model mentioned above in the Gaussian case, the problem
of identification of an exponent of self-similarity $E_0$ (which in some case is not unique) and of an Hurst
index $H_0$ is an open problem.

The first step in the resolution of this question involves an identification of some specific features of exponents and indices which can be recovered on sample paths. This paper is a first step : we will prove that from the regularity point of view these exponents and Hurst indices satisfy what we call optimality properties. More precisely, we prove that (see Theorem~\ref{ThOptim2}), the Hurst index $H_0$ maximizes the local critical exponent of the field in specific functional spaces related with the anisotropy matrix $E_0$ among all possible critical exponents in general anisotropic functional spaces.

Therefore, the results of the present paper open the way to the following strategy to recover the Hurst index. One first have to consider a discretized version of the set of all possible anisotropies. In each case an estimator of the critical exponent related with these anisotropies has to be given. Therefore, one has to locate the maximum of all these estimators--which can be based on anisotropic quadratic variations--and to identify the corresponding values of the anisotropy. Combining the resluts of this paper and these of \cite{38}, our approach can thus be turned in an effective algorithm for the estimation of the anisotropy of self similar textures (see~\cite{28,29}). The study of the estimators related to critical indices in anisotropic Besov spaces from a statistical point of view will be the purpose of a forthcoming paper.

Our optimality result comes from sample paths properties of the model under study in an
anisotropic setting. This approach is natural : In~\cite{19}, Kamont studied the regularity of the sample paths
of the well-known anisotropic Fractional Brownian Sheet in
anisotropic H\"{o}lder spaces related to Fractional Brownian Sheet. Moreover, some results of
regularity in specific anisotropic H\"{o}lder spaces related to matrix $E_0$
have already be established for operator scaling self-similar
random fields (which may be not Gaussian) in~\cite{5} or in the
more general setting of strongly non deterministic anisotropic
Gaussian fields in~\cite{37}. We then extend already existing results by measuring smoothness in general anisotropic spaces not necessarily
related to the exponent matrix $E_0$ of the field.

This paper is organized as follows. In Section~\ref{SecPres}, we briefly recall some facts about Operator Scaling Random Gaussian Fields (OSRGF) and describe the construction of~\cite{6} of the model. In Section \ref{SecAnisBesDef}, we present the different concepts used for measuring smoothness in an anisotropic setting and especially anisotropic Besov spaces. Section~\ref{SecOptim} is devoted to the statement of our optimality and regularity results. Finally, Section~\ref{SecProof} contains proofs of the results stated in Section~\ref{SecOptim}.

In the sequel, we will use some notations. For any matrix $M$
\[
\rho_{\min}(M)=\min\limits_{\lambda\in Sp(M)}(|\mathrm{Re}(\lambda)|),\;
\rho_{\max}(M)=\max\limits_{\lambda\in Sp(M)}(|\mathrm{Re}(\lambda)|)\;,
\]
where $\mathrm{Sp}(M)$ denotes the spectrum of matrix $M$.\\
For any real $a>0$, $a^{M}$ denotes the matrix
$$
a^{M}=\exp(M\log(a))=\sum\limits_{k\geq
0}\frac{M^{k}\log^{k}(a)}{k!}.
$$
In the following pages, we denote
$\mathcal{E}^{+}$ the collection of matrices of $M_d(\mathbb{R})$ whose
eigenvalues have positive real part.
\section{Presentation of the studied model}\label{SecPres}
The existence of operator scaling stable random fields, that is random fields satisfying relationship~(\ref{selfsimilar}), is proved in~\cite{6}. The following theorem (Theorem $4.1$ and Corollary $4.2$ of~\cite{6}) completes this result by yielding a practical way to construct a Operator Scaling Stable Random Field (OSRF) with stationary increments for any $E_0\in \mathcal{E}^{+}$ and $H_0\in (0,\rho_{\min}(E_0))$. We state it only in the Gaussian case, having in mind the problem of the estimation of the Hurst index $H_0$ and the anisotropy $E_0$.

\begin{theorem}\label{ThConstr}
Let $E_0$ be in $\mathcal{E}^{+}$ and $\rho$ a continuous function with positive values such that for all $x\neq 0$, $\rho(x)\neq 0$. Assume that $\rho$ is $\;E_0^{\!\!\!\!\!\!\!t}$--homogeneous, that is :
\[
\forall a>0,\,\forall
\xi\in\mathbb{R}^{d},\,\rho(a^{E_0^{\!\!\!\!\!\!\!t}}\xi)=a\rho(\xi)\;.
\]
Then the Gaussian field
\begin{equation}\label{EqDefX}
X_{\rho}(x)=\displaystyle\int_{\mathbb{R}^{d}}(\rme^{\rmi<x,\xi>}-1)\rho(\xi)^{-H_0-\frac{\mathrm{Tr}(E_0)}{2}}\rmd\widehat{W}(\xi)\;,
\end{equation}
exists and is stochastically continuous if and only if $H_0\in
(0,\rho_{\min}(E_0))$. Moreover this field has the following properties :
\begin{enumerate}
\item Stationary increments :
\[
\forall h\in\mathbb{R}^{d},\,\{X_{\rho}(x+h)-X_{\rho}(h)\}_{x\in\mathbb{R}^{d}}\overset{(fd)}{=}\{X_{\rho}(x)\}_{x\in\mathbb{R}^{d}}\;.
\]
\item The operator--scaling relation (\ref{selfsimilar}) is satisfied.
\end{enumerate}
\end{theorem}
\begin{remark}
The assumption of homogeneity on the function $\rho$ is necessary to recover linear self-similarity properties of the Gaussian field $\{X_{\rho}(x)\}_{x\in\mathbb{R}^d}$. The assumption of continuity on $\rho$ ensures that the constructed field is stochastically continuous.
\end{remark}
\begin{remark} In general, the couple $(H_0,E_0)$ of an OSRF is not unique. Indeed, if $H_0$ and $E_0$ are
respectively an Hurst index and an exponent of the OSRF $\{X(x)\}_{x\in\mathbb{R}^d}$, then
for any $\lambda>0$ so do $\lambda H_0$ and $\lambda E_0$.\\
Uniqueness of the Hurst index $H_0$ can be recovered by choosing a normalization for $E_0$, for example $\mathrm{Tr}(E_0)=d$.
However, even under this assumption, $E_0$ is not necessarily unique. Nevertheless remark that, under the assumption $\mathrm{Tr}(E_0)=d$, two anisotropies of an OSRF have necessarily the same real diagonalizable part (see Section~\ref{SecRegLoc} for a definition). We refer to Remark 2.10 of \cite{6}
for more details on the structure of the set of exponents of an
OSRF.
\end{remark}

Remark that Theorem~\ref{ThConstr} relies on the existence of $\;E_0^{\!\!\!\!\!\!\!\!t}$--homogeneous functions. Constructions of such functions have been proposed in~\cite{6} via an integral formula (Theorem~2.11). An alternative construction, more fitted for numerical simulations, can be found in~\cite{11}.

\section{Anisotropic concepts of smoothness}\label{SecAnisBesDef}
Our main goal here is to study the sample paths properties of this class of
Gaussian fields in suitable anisotropic functional spaces. This
approach is quite natural (see~\cite{19,5}) since the studied model is anisotropic. To this end, suitable concepts of anisotropic smoothness are needed. The aim of this section is to give some background about the appropriate anisotropic functional spaces : Anisotropic
Besov spaces. These spaces generalize classical (isotropic) Besov
spaces and have been studied in parallel with them (see~\cite{8,9} for a complete account on
the results presented in this section). The definition of anisotropic Besov spaces is based on the concept of pseudo-norm.
We first recall some well--known facts about pseudo-norms which can be found with more details in~\cite{21}.
\subsection{Preliminary results about pseudo-norms}\label{SecThree}
In order to introduce anisotropic functional spaces, an
anisotropic topology on $\mathbb{R}^d$ is needed. We need to
introduce a slight variant of the notion of pseudo--norm
introduced in~\cite{21}, fitted to the case of discrete
dilatations.
\begin{definition}
Let $E\in \mathcal{E}^+$. A function $\rho$ defined on $\mathbb{R}^d$ is a $(\mathbb{R}^d,E)$ pseudo-norm if
it satisfies the three following properties :
\begin{enumerate}
\item $\rho$ is continuous on $\mathbb{R}^d$,
\item $\rho$ is $E$-homogeneous, {\it i.e.}
$\rho(a^E x)= a\rho(x) \quad \forall x \in \mathbb{R}^d, \, \forall a>0$,
\item $\rho$ is strictly positive on $\mathbb{R}^d \setminus \{0\}$ .
\end{enumerate}
\end{definition}

For any $(\mathbb{R}^d,E)$ pseudo--norm, define the anisotropic sphere $S_0^E(\rho)$ as
\begin{equation}
S_0^E(\rho)= \{x \in \mathbb{R}^d; \, \rho(x)=1 \}\;.
\end{equation}
\begin{figure}[H]\label{Fig2}
  \begin{minipage}[c]{.32\linewidth}
    \centering
 \includegraphics[angle=0,width=.88\textwidth]{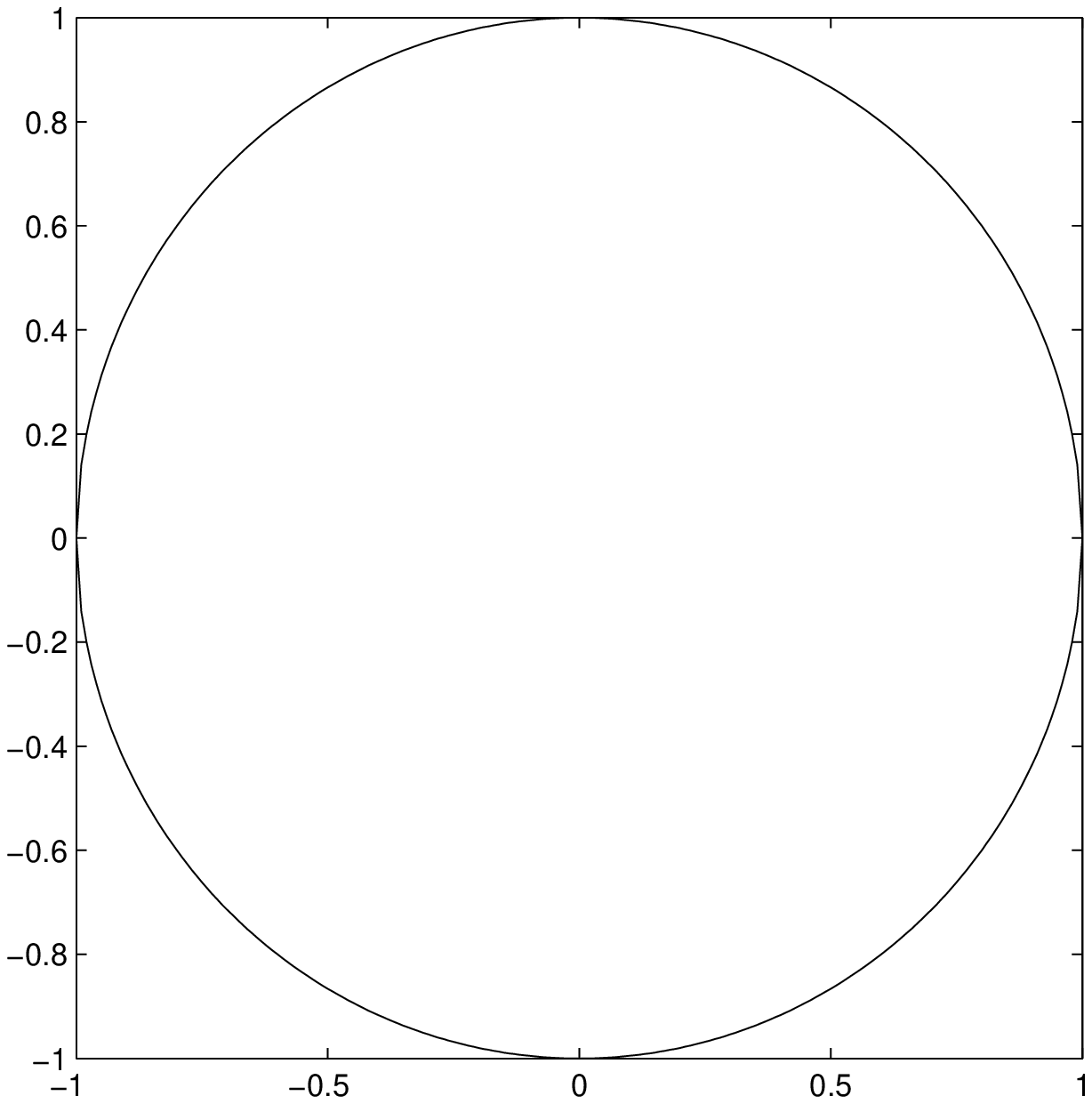}
\begin{center}$E=\begin{pmatrix}1&0\\0&1\end{pmatrix}$\end{center}
\end{minipage}\hfill
  \begin{minipage}[c]{.32\linewidth}
    \hspace{.025\textwidth}
      \includegraphics[angle=0,width=.95\textwidth]{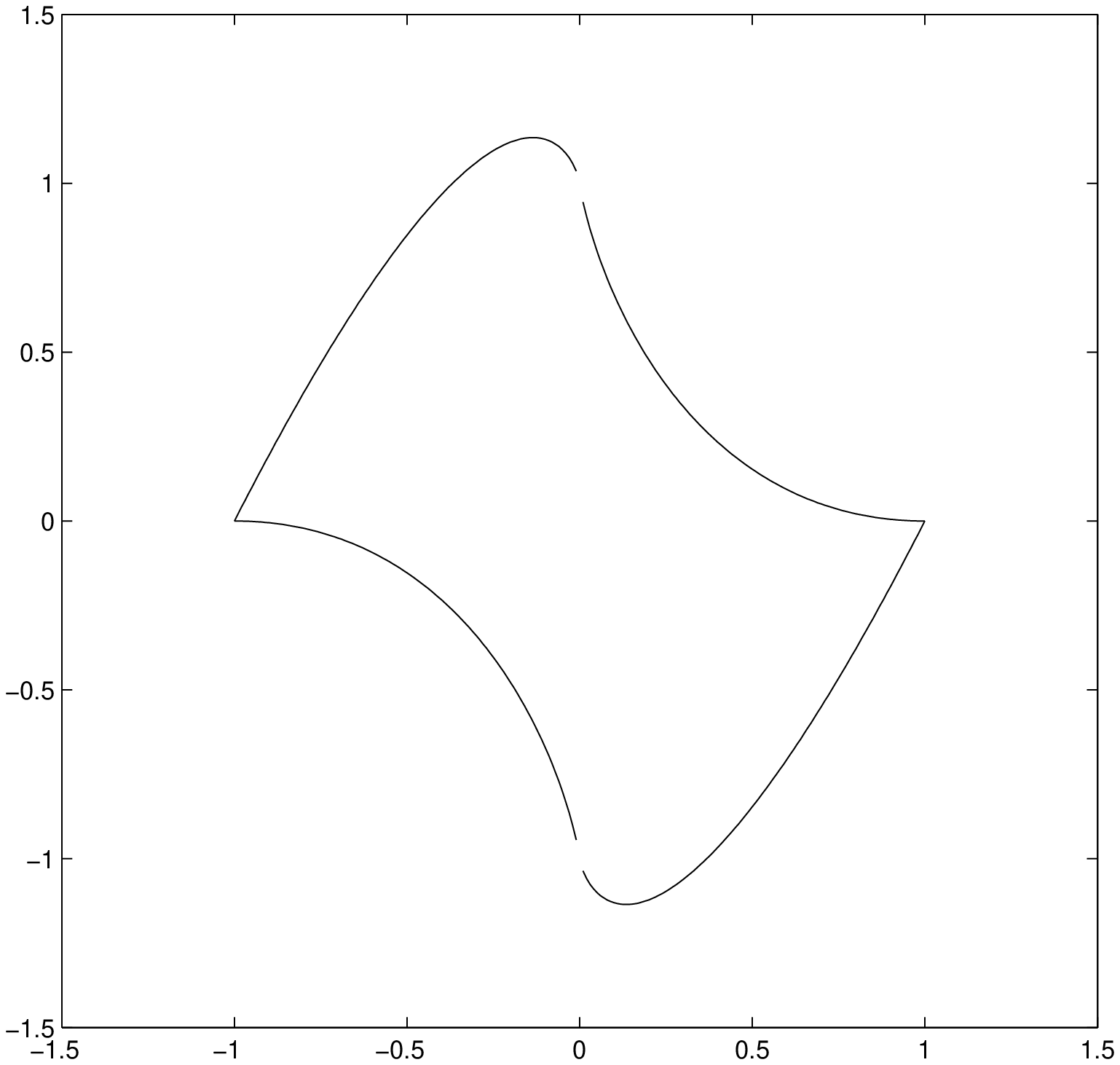}
      \begin{center}$E=\begin{pmatrix}1&1\\0&1\end{pmatrix}$\end{center}
\end{minipage}
 \begin{minipage}[c]{.32\linewidth}
    \hspace{.025\textwidth}
       \includegraphics[angle=0,width=.97\textwidth]{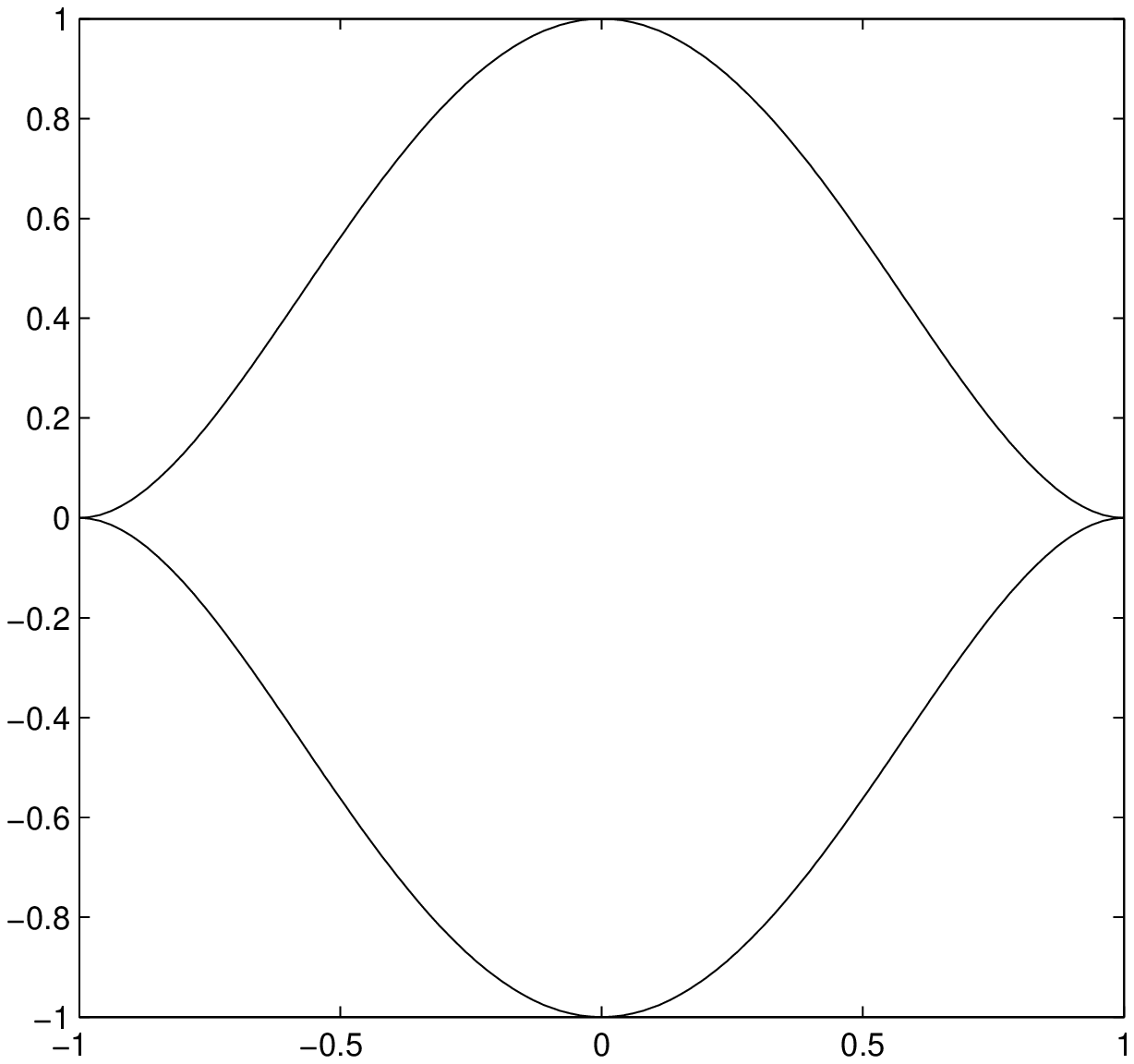}
       \begin{center}$E=\begin{pmatrix}1&0\\0&4\end{pmatrix}$\end{center}
\end{minipage}
\begin{center}\begin{caption}\; Examples of anisotropic spheres for different anisotropies.\end{caption}\end{center}
\end{figure}

\begin{proposition}
For all $x \in \mathbb{R}^d \setminus \{0\} $, there exists an
unique couple $(r, \theta) \in \mathbb{R}^*_+ \times
S_0^{E}(\rho)$ such that $x = r^{E} \theta$.

Moreover $S_0^{E}(\rho)$ is a compact of $\mathbb{R}^d$ and the map
$$
(r, \theta) \to x =r^{E} \theta\;,
$$
is an homeomorphism from $\mathbb{R}^*_+ \times S_0^{E}(\rho)$ to
$\mathbb{R}^d \setminus \{0\}$.
\end{proposition}
The term ``pseudo-norm'' is justified by the following proposition :
\begin{proposition}
Let $\rho$ a $(\mathbb{R}^d,E)$ pseudo-norm. There exists a constant $C>0$ such that
\begin{equation}\label{e:triangineq}
\rho(x+y) \le C (\rho(x)+\rho(y)), \quad \forall x, \, y \in
\mathbb{R}^d\;.
\end{equation}
\end{proposition}
The following key property allows to define an anisotropic topology on $\mathbb{R}^d$ based on pseudo-norms and then anisotropic functional spaces :
\begin{proposition}\label{eqpseudonormes} Let $\rho_1$ and $\rho_2$
be two $(\mathbb{R}^d, E)$ pseudo-norms. They are equivalent in the following sense : There exists a constant $C>0$ such that
$$
\frac{1}{C} \rho_1(x) \le \rho_2(x) \le C\rho_1(x), \quad
\forall x \in \mathbb{R}^d\;.
$$
In particular, two different $(\mathbb{R}^d, E)$ pseudo-norms define the same topology on $\mathbb{R}^{d}$.
\end{proposition}

\subsection{Anisotropic Besov spaces}
Let $E\in\mathcal{E}^+$ and $| \cdot |_{\;E^{\!\!\!\!\!\!\!t}}\;\;$ a fixed
$(\mathbb{R}^d, E^{\!\!\!\!\!\!\! t}\;\;)$--pseudo-norm. For $x_0 \in \mathbb{R}^d$ and
$r>0$, $B_{\vert \cdot \vert_{\;E^{\!\!\!\!\!\!t}}}\;(x_0,r)$ denotes the
anisotropic ball of center $x_0$ and radius $r$, namely
$$
B_{|\cdot|_{\;E^{\!\!\!\!\!\!t}}}\;(x_{0},r)=\{x\in\mathbb{R}^{d},|x-x_{0}|_{\;E^{\!\!\!\!\!\!\!t}}\;\leq r\}\;.
$$
\begin{definition}
Let $\psi_{0}^{E}\in\mathcal{S}(\mathbb{R}^{d})$ be such that
\[
\left\{
\begin{array}{l}
\widehat{\psi_{0}^{E}}(\xi)=1\mbox{ if }|\xi|_{\,E^{\!\!\!\!\!\!t}}\;\leq
1,\\
\widehat{\psi_{0}^{E}}(\xi)=0\mbox{ if }|\xi|_{\,E^{\!\!\!\!\!\!t}}\;\geq 2\;.
\end{array}
\right.
\]
For any positive integer $j$, set
\[
\widehat{\psi_{j}^{E}}(\xi)=\widehat{\psi_{0}^{E}}(2^{-j\,E^{\!\!\!\!\!\!t}}\;\xi)-\widehat{\psi_{0}^{E}}(2^{-(j-1)\,E^{\!\!\!\!\!\!t}}\;\xi)\;.
\]
Then
\[
\sum\limits_{j=0}^{+\infty}\widehat{\psi_{j}^{E}}\equiv 1\;,
\]
is an anisotropic partition of the unity satisfying
$\mathrm{supp}(\widehat{\psi_{j}^{E}})\subset
B_{|\cdot|_{\;E^{\!\!\!\!\!\!t}}}\;(0,2^{j+1})\setminus
B_{|\cdot|_{\;E^{\!\!\!\!\!\!t}}}\;(0,2^{j-1})$.
\end{definition}
The anisotropic Besov spaces
$B^{s}_{p,q}(\mathbb{R}^{d},E)$ are then defined as follows~:
\begin{definition}Let $0<p,q\leq\infty$ and  $s\in\mathbb{R}$. Define
\begin{equation}\label{e:defnorm1}
\|f\|_{B^{s}_{p,q}(\mathbb{R}^{d},E)}=
\left(\sum\limits_{j=0}^{\infty}2^{jsq}\|f*\psi_{j}^{E}\|^{q}_{L^{p}(\mathbb{R}^{d})}\right)^{1/q}\;.
\end{equation}
Then
\[
B^{s}_{p,q}(\mathbb{R}^{d},E)=\{f\in\mathcal{S}'(\mathbb{R}^{d}),\|f\|_{B^{s}_{p,q}(\mathbb{R}^{d},E)}<+\infty\}\;.
\]
The matrix $E$ is called the anisotropy of the Besov space
$B^{s}_{p,q}(\mathbb{R}^{d},E)$.\\
In a more general way, if $\beta\in\mathbb{R}$, define
\[
\|f\|_{B^{s}_{p,q,|\log|^{\beta}}(\mathbb{R}^{d},E)}=
\left(\sum\limits_{j=0}^{\infty}j^{-\beta q}2^{jsq}\|f*\psi_{j}^{E}\|^{q}_{L^{p}(\mathbb{R}^{d})}\right)^{1/q}\;.
\]
Then
\[
B^{s}_{p,q,|\log|^{\beta}}(\mathbb{R}^{d},E)=
\{f\in\mathcal{S}'(\mathbb{R}^{d}),\|f\|_{B^{s}_{p,q,|\log|^{\beta}}(\mathbb{R}^{d},E)}<+\infty\}\;.
\]
\end{definition}
\begin{remark}
One can prove that this definition is independent of the choice of the function $\psi_{0}^{E}$ involved in the definition of the Besov space $B^{s}_{p,q}(\mathbb{R}^{d},E)$.
\end{remark}
\begin{remark}\label{rem:comp} Let $E\in\mathcal{E}^+$ and $|\cdot|_{\;E^{\!\!\!\!\!\!t}}\;\,$ a $(\mathbb{R}^d, \;E^{\!\!\!\!\!\!\!t}\;\;)$ pseudo-norm. For any
$\lambda>0$, $|\cdot|^{1/\lambda}_{\;E^{\!\!\!\!\!\!\!\!t}}\;$ is a
$(\mathbb{R}^d, \lambda \;E^{\!\!\!\!\!\!\!\!t}\;\;)$ pseudo-norm. Hence for any $s>0$,
$B^{\lambda s}_{p,q}(\mathbb{R}^d,\lambda E)=B^{s}_{p,q}(\mathbb{R}^d, E)$.
\end{remark}

So, without loss of generality, we assume in the sequel that
$\mathrm{Tr}(E)=d$. We then define
\[
\mathcal{E}^+_d=\{E\in\mathcal{E}^+,\,\mathrm{Tr}(E)=d\}\;,
\]
where $\mathcal{E}^+$ is the collection of $d\times d$ matrices whose eigenvalues have positive real part.
As it is the case for isotropic spaces, anisotropic H\"older
spaces $\mathcal{C}^s(\mathbb{R}^d, E)$ can be defined as particular anisotropic Besov
spaces.
\begin{definition}
Let $s$ be in $\mathbb{R}$ and $\beta\in\mathbb{R}$. The anisotropic
H\"older spaces ${\mathcal{C}}^s(\mathbb{R}^d, E)$ and
${\mathcal{C}}^s_{\vert \log \vert^N}(\mathbb{R}^d, E)$ are defined by
\[
\mathcal{C}^{s}(\mathbb{R}^{d},E)=B^{s}_{\infty,\infty}(\mathbb{R}^{d},E)
\quad {\text{and}} \quad
\mathcal{C}^{s}_{|\log|^{\beta}}(\mathbb{R}^{d},E)=B^{s}_{\infty,\infty,|\log|^{\beta}}(\mathbb{R}^{d},E)\;.
\]
\end{definition}

\begin{proposition}\label{PropCarDifHold}
Let $0<s<\rho_{\min}(E)$ and $\beta\in\mathbb{R}$. Then the two norms
\begin{itemize}
\item
\[
\|f\|_{L^{\infty}(\mathbb{R}^{d})}+\sup_{|h|_{E}\leq
1}\sup_{x\in\mathbb{R}^{d}}\left(\frac{|f(x+h)-f(x)|}{|h|_{E}^{s}|\log(|h|_E)|^\beta}\right)\;,
\]
\item $\|f\|_{B^{s}_{\infty,\infty,|\log|^\beta}}$ defined by~(\ref{e:defnorm1}),
\end{itemize}
are equivalent in $\mathcal{C}^{s}_{|\log|^{\beta}}(\mathbb{R}^{d},E)$.\\
\end{proposition}
\begin{remark} Anisotropic H\"{o}lder spaces admit a
characterization by finite
differences of order $M\geq 1$ under the general assumption $s>0$. Here, we only need to deal with the case $0<s<\rho_{\min}(E)$ and have thus stated Proposition \ref{PropCarDifHold} in this special setting.\\
\end{remark}
Let us comment Proposition~\ref{PropCarDifHold}. Let
$0<s<\rho_{\min}(E)$ and $N\in\mathbb{R}$. A bounded function $f$
belongs to $\mathcal{C}^{s}_{|\log|^{\beta}}(\mathbb{R}^{d},E)$ if and
only if for any $r\in (0,1)$, $\Theta\in S_0^{E}(|\cdot|_E)$ and
$x\in \mathbb{R}^d$, one has
\[
|f(x+r^E \Theta)-f(x)|\leq C_0 r^s|\log(r)|^\beta\;,
\]
for some $C_0>0$.\\
Hence, a function $f$ belongs to the H\"{o}lder space
$\mathcal{C}^{s}_{|\log|^{\beta}}(\mathbb{R}^{d},E)$ if and only if its
restriction $f_{\Theta}$ along any parametric curve of the form
\[
r>0\mapsto r^E \Theta\;,
\]
with $\Theta\in S_0^{E}(|\cdot|_E)$ is in the usual H\"{o}lder
space $\mathcal{C}^{s}_{\vert \log \vert^\beta}(\mathbb{R})$ and
$\|f_{\Theta}\|_{\mathcal{C}^{s}_{\vert \log \vert^\beta}(\mathbb{R})}$
does not depend on $\Theta$. Roughly speaking, the anisotropic
``directional'' regularity in any anisotropic ``direction'' has to
be larger than $s$. In other words, we replace straight lines of
isotropic setting by curves with parametric equation $r>0\mapsto r^E \Theta$ adapted to anisotropic setting.\\
To state our optimality results we need a local version of
anisotropic Besov spaces~:
\begin{definition}
Let $E \in \mathcal{E}^{+}$ be a fixed anisotropy, $0<p,q\leq \infty$, $0<s<\infty$ and $f\in L^{p}_{loc}(\mathbb{R}^{d})$.\\
The function $f$ belongs to
$B^{\alpha}_{p,q,loc}(\mathbb{R}^{d},E)$ if for any $\varphi\in
\mathcal{D}(\mathbb{R}^{d})$, the
function $\varphi f$ belongs to $B^{\alpha}_{p,q}(\mathbb{R}^{d},E)$.\\
The spaces $B^{\alpha}_{p,q,|\log|^\beta,loc}(\mathbb{R}^{d},E)$ can be defined in an analogous way for any
$0<p,q\leq \infty$, $0<s<\infty$, $\beta\in\mathbb{R}$.\\

The anisotropic local critical exponent in anisotropic Besov
spaces $B^{s}_{p,q}(\mathbb{R}^d,E)$ of $f\in
L^p_{loc}(\mathbb{R}^{d})$ is then defined by
\[
\alpha_{f,loc}(E,p,q)=\sup\{s,\,f\in
B^{s}_{p,q,loc}(\mathbb{R}^d,E)\}\;.
\]
\end{definition}
In the special case $p = q = \infty$, this exponent is also called
the anisotropic local critical exponent in anisotropic H\"{o}lder
spaces of $f\in L^\infty_{loc}(\mathbb{R}^{d})$ and is denoted by
$\alpha_{f,loc}(E)$.

\section{Statement of our results}\label{SecOptim}
In what follows, we are given $E_0\in\mathcal{E}^+$ and $\rho_{E_0}$
a $(\mathbb{R}^d,\; E_0^{\!\!\!\!\!\!\! t}\;)$ pseudo-norm. We denote
$\{X_{\rho_{E_0},H_0}(x)\}_{x\in\mathbb{R}^d}$ the OSRGF with
exponent $E_0$ and Hurst index $H_0$ defined by~(\ref{EqDefX})
with $\rho=\rho_{E_0}$.

We first state our optimality result and
characterize in some sense an anisotropy $E_0$ and an Hurst index
of the field $\{X_{\rho_{E_0},H_0}(x)\}_{x\in\mathbb{R}^d}$. These
results come from an accurate study of sample paths properties of
the OSRGF $\{X_{\rho_{E_0},H_0}(x)\}_{x\in\mathbb{R}^d}$ in
anisotropic Besov spaces (see in Section~\ref{SecSample}).

We assume --without loss of generality --that
$E_0\in\mathcal{E}^+_d$, namely that all the eigenvalues of $E_0$ have a positive real part and that $\mathrm{Tr}(E_0)=d$. Our results will be based on a comparison
between the topology related to the pseudo--norm $\rho_{E_0}$
involved in the construction of the Gaussian field
$\{X_{\rho_{E_0},H_0}(x)\}_{x\in\mathbb{R}^d}$ and this of the
analyzing spaces $B^s_{p,q}(\mathbb{R}^d,E)$. To be able to
compare these two topologies, we also assume that
$E\in\mathcal{E}^+_d$.

The main result of this paper is the following one~:
\begin{theorem}\label{ThOptim2}
Let $(p,q)\in [1,+\infty]^2$ and $E_0\in\mathcal{E}^+_d$. Then
almost surely
\begin{eqnarray*}
\alpha_{X_{\rho_{E_0},H_0},loc}(E_0,p,q)&=&\sup\{\alpha_{X_{\rho_{E_0},H_0},loc}(E,p,q),E\in
\mathcal{E}^{+}_d, E\mbox{ commuting with }E_0\}\\
&=&H_0\;,
\end{eqnarray*}
that is the value $E=E_0$ maximizes the anisotropic local critical
exponent of the OSRGF $\{X_{\rho_{E_0},H_0}(x)\}_{x\in\mathbb{R}^d}$ among
all possible anisotropic local critical exponent in anisotropic
Besov spaces with an anisotropy $E$ commuting with $E_0$.
\end{theorem}
\begin{remark} Since $E$ and $E_0$ are commuting, these matrices admit the same
spectral decomposition. Hence, in fact we proved that any
anisotropy $E_0$ maximize the critical exponent among
matrices having the same spectral decomposition. Thus, in the general case, we implicitly assumed
that the spectral decomposition of anisotropy matrix is known.
In dimension two, we have a stronger optimality result about anisotropy
$E_0$ and Hurst index $H_0$, involving matrices of $\mathcal{E}^{+}_d$ which do not commute necessilary.
\end{remark}
To prove Theorem~\ref{ThOptim2}, we investigate the local regularity
of the sample paths of $\{X_{\rho_{E_0},H_0}(x)\}_{x\in\mathbb{R}^d}$ in general anisotropic Besov spaces. But before any statement, we first need some background about the concept of real diagonalizable part of a square matrix. This notion is based on real additive Jordan decomposition of a square matrix (see for e.g. to Lemma 7.1 chap 9 of \cite{17} where a multiplicative version of Proposition~\ref{PropDecJAdd} is given) :
\begin{proposition}\label{PropDecJAdd}
Any matrix $M$ of
$M_{d}(\mathbb{R})$ can be decomposed into a sum of three commuting real matrices
\[
M=D+S+N\;,
\]
where $D$ is a diagonalizable matrix in $M_{d}(\mathbb{R})$, $S$ is a diagonalizable matrix in $M_{d}(\mathbb{C})$ with zero or
imaginary complex eigenvalues, and $N$ is a nilpotent matrix. Matrix $D$ is called the real diagonalizable part of
$M$, $S$ its imaginary semi-simple part, and $N$ its
nilpotent part.
\end{proposition}

Now we are given two {\bf commuting} matrices $E_0$, $E$
of $\mathcal{E}^{+}_{d}$. Let $D_0$ (resp $D$) be the real diagonalizable part of matrix $E_0$ (resp $E$). Since matrices $E_0$ and $E$ are commuting, so do matrices $D_0$
and $D$. Furthermore, matrices $D_0$ and $D$ are diagonalizable in
$M_{d}(\mathbb{R})$ then they are simultaneously diagonalizable.
Up to a change of basis, we may assume that
$D_0$ and $D$ are two diagonal matrices. More precisely, suppose that
\begin{equation}\label{EqMatricesDDp}
D_0=\begin{pmatrix}\lambda_{1}^{0}Id_{d_{1}}&
&0\\&\ddots&\\0&&\lambda_{m}^{0}Id_{d_{m}}\end{pmatrix},\,D=\begin{pmatrix}\lambda_{1}Id_{d_{1}}&
&0\\&\ddots&\\0&&\lambda_{m}Id_{d_{m}}\end{pmatrix},
\end{equation}
with
\begin{equation}\label{EqValeursPropres}
\frac{\lambda_{m}}{\lambda_{m}^{0}}\leq\cdots\leq
\frac{\lambda_{1}}{\lambda_{1}^{0}}.
\end{equation}
Since $\mathrm{Tr}(E_0)=\mathrm{Tr}(E)=d$, one has $\lambda_{m}/\lambda_{m}^{0}\leq 1$. \\

The regularity results about sample path of the field
$\{X_{\rho_{E_0},H_0}(x)\}_{x\in\mathbb{R}^d}$ are summed up in the following theorem.
\begin{theorem}\label{PropRegBesovAutreAnisotr}Let $1\leq p\leq +\infty$, $1\leq q\leq +\infty$.
Almost surely the anisotropic local critical exponent $\alpha_{X_{\rho_{E_0},H_0},loc}(E,p,q)$ in anisotropic Besov spaces
$B^{s}_{p,q}(\mathbb{R}^d,E)$ of the OSRGF
$\{X_{\rho_{E_0},H_0}(x)\}_{x\in\mathbb{R}^d}$ satisfies
\[
\alpha_{X_{\rho_{E_0},H_0},loc}(E,p,q)=\frac{\lambda_{m}H_0}{\lambda_{m}^{0}}\leq H_0\;.
\]
In particular, in the special case $E=E_0$, one has $\alpha_{X_{\rho_{E_0},H_0},loc}(E,p,q)=H_0$.

\end{theorem}
In other words Theorem~\ref{PropRegBesovAutreAnisotr} asserts that
when one measures local regularity of the sample paths along
anisotropic directions different from those associated to an
anisotropy of the field $E_0$, one loses smoothness. The further the anisotropic
direction of measure from the genuine anisotropic direction
associated to the field are, the smaller the anisotropic local critical exponent is.
This anisotropic local critical exponent can take any value in the range $(0,H_0]$.\\

The special case $p=q=+\infty$ yields us
the following result about anisotropic H\"{o}lderian regularity of the
sample paths.
\begin{corollary}\label{PropCompHoldAnis2}Almost surely the anisotropic local critical exponent of the sample paths of $\{X_{\rho_{E_0},H_0}(x)\}_{x\in\mathbb{R}^d}$
in anisotropic H\"{o}lder spaces equals
$(\lambda_{m} H_0)/\lambda_{m}^{0}$ and is always lower than $H_0$. In particular, if $E=E_0$ this critical exponent equals the Hurst index $H_0$.
\end{corollary}
\begin{remark}
This estimate on anisotropic local critical exponent was already known in the case $E=E_0$ (see~\cite{5}).
\end{remark}

Theorem \ref{PropRegBesovAutreAnisotr} allows us to obtain
regularity results which extend those proved in the case
$p=q=\infty$ in the usual isotropic setting. Since matrices $E_0$
and $Id$ are commuting, we can apply the above result to the case
$E=Id$. Note that in this case $\lambda_{m}^{0}=\rho_{\max}(E_0)$.
We obtain the following proposition :
\begin{proposition}\label{PropCasIso}
Almost surely the local critical exponent of the sample paths of
$\{X_{\rho_{E_0},H_0}(x)\}_{x\in\mathbb{R}^d}$ in classical Besov
spaces equals $H_0/\rho_{\max}(E_0)$.

In particular, for $p=q=\infty$, almost surely the local critical exponent of the sample paths of $\{X_{\rho_{E_0},H_0}(x)\}_{x\in\mathbb{R}^d}$ in classical H\"{o}lder spaces
equals $H_0/\rho_{\max}(E_0)$.\\
\end{proposition}
\begin{remark} In the special case $p=q=\infty$, we recover already known results about classic H\"{o}lderian
regularity (see Theorem~$5.4$ of~\cite{6}). Recall that this
theorem is based on directional regularity results about the
Gaussian field $\{X_{\rho_{E_0},H_0}\}$ and comes from an estimate
of the variogram
$v_{X_{\rho_{E_0},H_0}}(h)=\mathbb{E}(|X_{\rho_{E_0},H_0}(h)|^2)$
along special directions related to the spectral decomposition of
matrix $E_0$. Here our approach is different and based on wavelet
technics.
\end{remark}
\section{Complements and proofs}\label{SecProof}
\subsection{Role of the real diagonalizable part of the
anisotropy $E$ of the analysing spaces $B^s_{p,q}(\mathbb{R}^d,E)$ }
We will first prove that measuring smoothness in the general Besov spaces $B^s_{p,q}(\mathbb{R}^d,E)$ may be deduced from the special case where the matrix $E$ is diagonalizable. To this end, we show the following embedding property:

\begin{proposition}\label{PropCompBesAnis1}
Assume that $E_1 \in {\mathcal{E}}^+_d$ and $E_2 \in {\mathcal{E}}^+_d$ have the same real diagonalizable part
$D$. Let $|\cdot|_{\;E_{1}^{\!\!\!\!\!\!\!t}}\;$ (resp $|\cdot|_{\;E_{2}^{\!\!\!\!\!\!\!t}}\;$)) a $(\mathbb{R}^{d},\;E_{1}^{\!\!\!\!\!\!\!t}\;)$ (resp $(\mathbb{R}^{d},\;E_{2}^{\!\!\!\!\!\!\!t}\;)$) pseudo--norm. Then for any $\alpha>0$ and any $(p,q)\in [1,+\infty]^2$ one has,
\begin{equation}\label{e:embed}
B^{\alpha}_{p,q,|\log|^{-\frac{d}{\rho_{\min}(D)}-1}}(\mathbb{R}^{d},E_{1})\hookrightarrow
B^{\alpha}_{p,q}(\mathbb{R}^{d},E_{2}) \hookrightarrow
B^{\alpha}_{p,q,|\log|^{\frac{d}{\rho_{\min}(D)}+1}}(\mathbb{R}^{d},E_{1})\;.
\end{equation}
\end{proposition}
As a direct consequence, we obtain Corollary~\ref{CorCompBesAnis1}.
\begin{corollary}\label{CorCompBesAnis1}
The anisotropic local critical exponent
\[
\alpha_{X,loc}(E,p,q)=\sup\{s>0,\,X(\cdot)\in B^s_{p,q,loc}(\mathbb{R}^d,E)\}\;,
\]
of any Gaussian field $\{X(x)\}_{x \in\mathbb{R}^d}$ in anisotropic Besov spaces
$B^{s}_{p,q}(\mathbb{R}^{d},E)$ depends only on the real
diagonalizable part of $E$.
\end{corollary}
Note that this result does not depend on the studied Gaussian
field but of the analyzing functional spaces. Hence, it does not
give any information about the anisotropic
properties of the field.

We now show Proposition~\ref{PropCompBesAnis1}. The proof of this result relies on the following lemma :
\begin{lemma}\label{LemComPNDiagCom}
Assume that $E_1$ and $E_2$ are two matrices of $\mathcal{E}^+_d$ having the same real diagonalizable part $D$.
Then there exists two positive constants $c_1$ and $c_2$ such
that, for all $x \in \mathbb{R}^d$,
\begin{equation}\label{e:ComPNDiagCom}
c_{1}|x|_{\;E_{2}^{\!\!\!\!\!\!\!t}}(1+|\log(|x|_{\;E_{2}^{\!\!\!\!\!\!\!t}})|)^{-\frac{d}{\rho_{\min}(D)}}\leq
|x|_{\;E_{1}^{\!\!\!\!\!\!\!t}}\leq
c_{2}|x|_{\;E_{2}^{\!\!\!\!\!\!\!t}}(1+|\log(|x|_{\;E_{2}^{\!\!\!\!\!\!\!t}})|)^{\frac{d}{\rho_{\min}(D)}}\;.
\end{equation}
\end{lemma}
{\bf Proof of Lemma~\ref{LemComPNDiagCom}.} Using polar coordinates associated to $\;E_{1}^{\!\!\!\!\!\!\!t}$, one has, for $x \in
\mathbb{R}^d$,
\[
x=r^{\;E_{1}^{\!\!\!\!\!\!\!t}}\Theta, (r,\Theta)\in \mathbb{R}^{*}_{+}\times
S_{0}^{\;E_{1}^{\!\!\!\!\!\!\!t}}\;(|\cdot|_{\;E_{1}^{\!\!\!\!\!\!\!t}}\;)\;.
\]
Denote $F_1=E_1-D$, $F_2=E_2-D$. Then
\begin{eqnarray*}
|x|_{\;E_{2}^{\!\!\!\!\!\!\!t}}& = & |r^{\;E_{2}^{\!\!\!\!\!\!\!t}}\cdot (r^{-D}r^{-\;F_2^{\!\!\!\!\!\!\!t}})\cdot (r^{D}r^{\;F_1^{\!\!\!\!\!\!\!t}}\Theta)|_{\;E_{2}^{\!\!\!\!\!\!\!t}}\\
& \leq& r|r^{-\;F_2^{\!\!\!\!\!\!\!t}}r^{\;F_1^{\!\!\!\!\!\!\!t}}\Theta|_{\;E_{2}^{\!\!\!\!\!\!\!t}}\;,
\end{eqnarray*}
because $F_1,F_2, D$ are pairwise commuting matrices.  Observe now that $F_1,F_2$ have only pure imaginary eigenvalues. Hence, by Lemma~$2.1$ of~\cite{6}, one deduces that for any $\varepsilon>0$
\begin{eqnarray*}
|x|_{\;E_{2}^{\!\!\!\!\!\!\!t}}&  \leq& C
r\max(|r^{-\;F_2^{\!\!\!\!\!\!t}}r^{\;F_1^{\!\!\!\!\!\!t}}\Theta|^{\frac{1}{\rho_{\min}(D)-\varepsilon}},
|r^{-\;F_2^{\!\!\!\!\!\!t}}r^{\;F_1^{\!\!\!\!\!\!t}}\Theta|^{\frac{1}{\rho_{\max}(D)+\varepsilon}})\;,
\end{eqnarray*}
where $|\cdot|$ denotes the usual Euclidean norm. Denote $\|\cdot\|$ an operator norm on $M_d(\mathbb{R})$. Since $\Theta$ belongs to the anisotropic sphere $S_{0}^{\;E_{1}^{\!\!\!\!\!\!\!t}}\;(|\cdot|_{\;E_{1}^{\!\!\!\!\!\!\!t}}\;)$ which is compact, one has
\begin{eqnarray*}
|x|_{\;E_{2}^{\!\!\!\!\!\!\!t}} & \leq&
Cr\max(\|r^{-\;F_{2}^{\!\!\!\!\!\!\!t}}r^{\;F_{1}^{\!\!\!\!\!\!\!t}}\|^{\frac{1}{\rho_{\min}(D)-\varepsilon}},
\|r^{-\;F_{2}^{\!\!\!\!\!\!\!t}}r^{\;F_{1}^{\!\!\!\!\!\!\!t}}\|^{\frac{1}{\rho_{\max}(D)+\varepsilon}})\\
 & \leq& Cr\;(1+|\log(r)|)^{\frac{d-1}{\rho_{\min}(D)-\varepsilon}}\\
 & \leq &Cr\;(1+|\log(r)|)^{\frac{d}{\rho_{\min}(D)}}\;,
\end{eqnarray*}
for $\varepsilon>0$ sufficiently small. We then proved Lemma~\ref{LemComPNDiagCom}. We now show Proposition~\ref{PropCompBesAnis1}.\\
{\bf Proof of Proposition~\ref{PropCompBesAnis1}.} Using two anisotropic Littlewood-Paley
analysis associated respectively to matrices $E_1,E_2$ and $D$ and the lemma above, we deduce~(\ref{e:embed}). Indeed, for any $i\in\{1,2\}$, let $(\psi_j^{E_i})_{j \in\mathbb{N}}$ an anisotropic Littlewood--Paley analysis of Besov spaces $B^{\alpha}_{p,q}(\mathbb{R}^d,E_i)$. By definition,
\[
\mathrm{supp}(\widehat{\psi_{1}^{E_{i}}})\subset\{\xi\in\mathbb{R}^{d},\,1\leq
|\xi|_{\;E_{i}^{\!\!\!\!\!\!\!t}}\;\leq 4\}\;,
\]
for $i\in\{1,2\}$. Then there exists some $j_0\in\mathbb{N}$ such that for any $j \in \mathbb{N}$, one has
\begin{eqnarray*}
\mathrm{supp}(\widehat{\psi_{j}^{E_{2}}})&\subset&\{\xi,2^{j-1}\leq
|\xi|_{\;E_{2}^{\!\!\!\!\!\!\!t}}\;\leq
2^{j+1}\}\\
&\subset&\bigcup\limits_{\ell=j-j_{0}-\frac{d\log_{2}(j)}{\rho_{\min}(D)}}^{j+j_{0}+\frac{d\log_{2}(j)}{\rho_{\min}(D)}}\{\xi,2^{\ell-1}\leq
|\xi|_{\;E_{1}^{\!\!\!\!\!\!\!t}}\leq 2^{\ell+1}\}\;.
\end{eqnarray*}
Hence
\[
\widehat{\psi_{j}^{E_{2}}}(\xi)\,\widehat{f}(\xi) =
\widehat{\psi}_{j}^{E_{2}}(\xi)\left(\sum\limits_{\ell=j-j_{0}-\frac{d\log_{2}(j)}{\rho_{\min}(D)}}^{j+j_{0}+\frac{d\log_{2}(j)}{\rho_{\min}(D)}}\widehat{\psi}_{\ell}^{E_{1}}(\xi)\,\widehat{f}(\xi)\right)\;.
\]
Define $q'$ the conjugate of $q$, that is the positive real satisfying $1/q+1/q'=1$. The last inequality and Cauchy--Schwartz inequality imply that for some $C>0$,
\begin{eqnarray*}
\|f*\psi_{j}^{E_{2}}\|_{L^{p}}^q&\leq&
C(\log_2 j)^{q/q'}\left(\sum\limits_{\ell=j-j_{0}-\frac{d\log_{2}(j)}{\rho_{\min}(D)}}^{j+j_{0}+\frac{d\log_{2}(j)}{\rho_{\min}(D)}}
\|\psi_{j}^{E_{2}}*(\psi_{\ell}^{E_{1}}*f)\|_{L^{p}}^q\right)\\
&\leq&C\;(\log_2 j)^{q/q'}\|\psi_0^{E_{2}}\|_{L^{1}}^q\left(\sum\limits_{\ell=j-j_{0}-\frac{d\log_{2}(j)}{\rho_{\min}(D)}}^{j+j_{0}+\frac{d\log_{2}(j)}{\rho_{\min}(D)}
} \|\psi_{\ell}^{E_{1}}*f\|_{L^{p}}^q\right)\;.
\end{eqnarray*}
Then we can give the following upper bound of
$\sum\limits_{j=1}^{J}2^{jsq}\|f*\psi_{j}^{E_{2}}\|_{L^{p}}^q$ :
\begin{eqnarray*}
\sum\limits_{j=1}^{J}2^{jsq}\|f*\psi_{j}^{E_{2}}\|_{L^{p}}^q&\leq&
C\sum\limits_{j=1}^{J}(\log_2 j)^{q/q'}2^{jsq}
\left(\sum\limits_{\ell=j-j_{0}-\frac{d\log_{2}(j)}{\rho_{\min}(D)}}^{j+j_{0}+\frac{d\log_{2}(j)}{\rho_{\min}(D)}}
\|(f*\psi_{\ell}^{E_{1}})\|_{L^{p}}^q\right)\\
&\leq&C\sum\limits_{\ell=1}^{J+j_{0}+\frac{d\log_{2}(J)}{\rho_{\min}(D)}}
\|f*\psi_{\ell}^{E_{1}}\|_{L^{p}}^q
\left(\sum\limits_{j=\ell-j_{0}-\frac{d\log_2(\ell)}{\rho_{\min}(D)}}^{\ell+j_{0}+\frac{d\log_2(\ell)}{\rho_{\min}(D)}}
(\log_2 j)^{q/q'}2^{jsq}\right)\\
&\leq&C\sum\limits_{\ell=1}^{J+j_{0}+\frac{d\log_{2}(J)}{\rho_{\min}(D)}}\|f*\psi_{\ell}^{E_{1}}\|_{L^{p}}^q 2^{\ell
sq}\ell^{d/\rho_{\min}(D)+1}\;.
\end{eqnarray*}
Let now $J$ tends to $\infty$. It yields the embedding
\[
B^{\alpha}_{p,q,|\log|^{-\frac{d}{\rho_{\min}(D)}-1}}(\mathbb{R}^{d},E_{1}) \hookrightarrow
B^{\alpha}_{p,q}(\mathbb{R}^{d},E_{2})\;.
\]
Permuting $E_1$ and $E_2$ yields the other inclusion.
\subsection{Local regularity in anisotropic Besov spaces of the studied
field}\label{SecRegLoc}
In the previous section, we proved that we
can restrict our study to diagonal Besov spaces. This point is
crucial for the proof of the regularity results stated in
Section~\ref{SecOptim}. Indeed it allows us to use tools that are
only defined in the diagonal case, as anisotropic multi-resolution
analysis and anisotropic wavelet bases. The aim of the following
subsection is to recall the constructions of these wavelet bases.
\subsubsection{Orthonormal Wavelet bases of (diagonal) anisotropic spaces}\label{SecBaseWAnis}

In this section, we assume that the anisotropy $D$ of the analyzing space is
diagonal (with positive eigenvalues), namely that
$$
D=\begin{pmatrix}\lambda_1&&0\\&\ddots&\\0&&\lambda_d\end{pmatrix}\;.
$$
In addition we also assume that $\mathrm{Tr}(D)=d$. Our main tool will be anisotropic multi--resolution analyses defined by Triebel
in~\cite{36}.

Let $\{V_j, j \geq 0\}$ be a one--dimensional multi--resolution analysis of
$L^2(\mathbb{R})$. Denote by $\psi^F$ (resp. $\psi^M$) the
corresponding scaling function (resp. wavelet function).

\begin{notation} We denote by $\{F,M\}^{d^{*}}$ the set
$$
\{F,M\}^{d^{*}}=\{F,M\}^{d}\setminus\{(F,\cdots,F)\}\;.
$$
For $j \in \mathbb{N}$, we define the set $I^j(D)$ of $\{F,M\}^d \times
\mathbb{N}^d$ in the following way.
\begin{itemize}
\item If $j=0$, $I^{0}(D)=\{((F,\cdots,F),(0,\cdots,0))\}$.\\
\item If $j\geq 1$, $I^{j}(D)$ is the set of all the elements
$(G,\gamma)$ with $G\in \{F,M\}^{d^{*}}$ and
$\gamma\in\mathbb{N}^{d}$ such that for any $r\in\{1,\cdots,d\}$ :
\[
\begin{array}{l}
\mbox{If }G_{r}=F,\gamma_{r}=[(j-1)\lambda_{r}],\\
\mbox{If }G_{r}=M,[(j-1)\lambda_{r}]\leq\gamma_{r}<[j\lambda_{r}]\;.
\end{array}
\]

Finally, for $j\in\mathbb{N}$ and $(G,\gamma)\in I^{j}(D)$, we
will denote by $D_{j,G,\gamma}$ the matrix defined by
\[
D_{j,G,\gamma}=\begin{pmatrix}\gamma_{1}&&0\\&\ddots&\\0&&\gamma_{d}\end{pmatrix}\;.
\]
\end{itemize}
\end{notation}

Finally, let us define the family of wavelets as follows. For
$j\in\mathbb{N}$, $(G,\gamma)\in I^{j}(D)$ and
$k\in\mathbb{Z}^{d}$, we set
\[
\Psi_{j,G,\gamma}^{k}(x)=(\psi^{(G)})(2^{D_{j,G,\gamma}}x-k)\;,
\]
with
\[
\psi^{(G)}=\psi_{G_{1}}\otimes\cdots\otimes\psi_{G_{d}}\;.
\]
The anisotropic wavelet bases yield a wavelet characterisation of
anisotropic Besov spaces (see~\cite{35} and~\cite{36},
Theorem $5.23$).
\begin{theorem}
\
\begin{enumerate}
\item The family
$\left\lbrace 2^{\frac{\mathrm{Tr}(D_{j,G,\gamma})}{2}}\Psi_{j,G,\gamma}^{k}, \, j\in\mathbb{N}, \, (G,\gamma)\in
I^{j}(D), \, k\in\mathbb{Z}^{d} \right\rbrace$ is an orthonormal basis of
$L^{2}(\mathbb{R}^{d})$. \item Let
$(\Psi^{j,G,\gamma}_{k})_{j\in\mathbb{N},(G,\gamma)\in
I^{j}(D),k\in\mathbb{Z}^{d}}$ be the family constructed from
$\psi_{F}$ and $\psi_{M}$ Daubechies wavelets with, for some $u
\in \mathbb{N}$,
\[
\psi_{F}\in\mathcal{C}^{u}(\mathbb{R}),\psi_{M}\in\mathcal{C}^{u}(\mathbb{R})\;.
\]
Let $0<p,q\leq\infty$ and $s,N\in\mathbb{R}$. There exists an
integer
 $u(s,p,D)$ such that if $u>u(s,p,D)$, for any tempered distribution $f$ the two following assertions are equivalent
 \begin{enumerate}
\item $f\in B^{s}_{p,q,|\log|^\beta}(\mathbb{R}^{d},D)$. \item $f=\sum
c_{j,G,\gamma}^{k}\Psi^{k}_{j,G,\gamma}$ with
\[
\sum\limits_{j,G,\gamma}j^{-\beta q} 2^{j(s-\frac{d}{p})q}
\left(\sum\limits_{k}|c_{j,G,\gamma}^{k}|^{p}\right)^{\frac{q}{p}}<+\infty\;,
\]
the convergence being in $\mathcal{S}'(\mathbb{R}^{d})$.
\end{enumerate}
The above expansion is then unique and
\begin{equation}\label{e:AnisotropicWC}
c_{j,G,\gamma}^{k}=<f,2^{\mathrm{Tr}(D_{j,G,\gamma})}\Psi^{k}_{j,G,\gamma}>\;.
\end{equation}
\end{enumerate}
\end{theorem}
\begin{remark}An analogous result is stated (see~\cite{36},
Theorem $5.24$) replacing Daubechies wavelets by Meyer wavelets.
In that case, $u=+\infty$.
\end{remark}
We now prove our regularity results about the sample path of
$\{X_{\rho_{E_0},H_0}(x)\}_{x\in\mathbb{R}^d}$ based on wavelet
characterization of Besov spaces.

\subsubsection{Local regularity of the field $\{X_{E_0,H_0}(x)\}_{x\in\mathbb{R}^{d}}$ in anisotropic
Besov spaces $B_{p,q}^{s}(\mathbb{R}^d,D_0)$} Assume that we are
given a Gaussian field $\{X_{E_0,H_0}(x)\}_{x\in\mathbb{R}^{d}}$
of the form~(\ref{EqDefX}) where $E_0\in\mathcal{E}^+_d$  and
$H_0\in (0,\rho_{\min}(E_0))$. The aim of this section is to prove
:
\begin{proposition}\label{PropLocRegAnisEgale} Let $1\leq p,q\leq+\infty$. Define $\delta$ on $(0,+\infty]$ as
follows :
\[
\delta(p)=\left\{
\begin{array}{l}
3/2\mbox{  if  }p=+\infty,\\
1\mbox{ otherwise}.
\end{array}
\right.
\]
Then one has
\begin{enumerate}
\item For any $\beta>1/q+d/\rho_{\min}(E_0)+\delta(p)$, almost surely,  the sample path of $\{X_{\rho_{E_0},H_0}(x)\}_{x\in\mathbb{R}^d}$
belongs to
$B^{H_0}_{p,q,|\log|^{\beta},loc}(\mathbb{R}^{d},D_0)$,
\item For $\beta=1/q+d/\rho_{\min}(E_0)+\delta(p)$, almost surely, the sample path of $\{X_{\rho_{E_0},H_0}(x)\}_{x\in\mathbb{R}^d}$ does
not belong to $B_{p,q,|\log|^{-\beta},loc}^{H_0}(\mathbb{R}^{d},D_0)$.
\end{enumerate}
\end{proposition}
Adapting to our setting a result of~\cite{21}, we first remark
that there exists
$\mathcal{C}^{\infty}(\mathbb{R}^{d}\setminus\{0\})$
$(\mathbb{R}^d,E_0)$ pseudo--norms
\begin{lemma}\label{lem:cinfty}Let $E_0\in\mathcal{E}^+$ and $\varphi$ be a $\mathcal{C}^{\infty}$ function
compactly supported in $\mathbb{R}^d\setminus\{0\}$. The function
$\rho$ defined on $\mathbb{R}^d$, by
\[
\rho(x)=\int_{\mathbb{R}^{d}}\varphi(a^{-E_0}x)\rmd a\;,
\]
is a $(\mathbb{R}^{d},E_0)$ pseudo-norm belonging to
$\mathcal{C}^{\infty}(\mathbb{R}^{d}\setminus \{0\})$.
\end{lemma}
We now prove that the belongness of the sample paths to anisotropic
Besov spaces of any OSRGF of the form
$\{X_{\rho,H_0}\}_{x\in\mathbb{R}^d}$ do not depend on the
$(\mathbb{R}^d,\;E_0^{\;\!\!\!\!\!\!\!\!\!t}\;)$ pseudo-norm $\rho$
involved in the construction of the field.
\begin{lemma}\label{lem:eqPN}
Let $E_0\in\mathcal{E}_d^+$ and $\rho_1,\rho_2$ two
$(\;E_0^{\!\!\!\!\!\!t},\mathbb{R}^d)$ pseudo--norms. Denote respectively
$\{X_1(x)\}_{x\in\mathbb{R}^{d}}$ and
$\{X_1(x)\}_{x\in\mathbb{R}^{d}}$ the two OSSRGF defined from
$\rho_1$ and $\rho_2$. Then, for any $s>0$, $\beta\in\mathbb{R}$, $(p,q)\in (0,\infty]^2$, a.s. $X_1$ belongs to $B^s_{p,q,|\log|^\beta}(\mathbb{R}^d,D_0)$ iff a.s.
$X_2$ belongs to $B^s_{p,q,|\log|^\beta}(\mathbb{R}^d,D_0)$.
\end{lemma}
{\bf Proof.} Remark first that using the same approach than in Lemma~2
of \cite{7} and an anisotropic version of Kolmogorov Centsov
Theorem we can prove that a.s.
\[
x\mapsto\int_{|\xi|\leq
R}(\mathrm{e}^{\mathrm{i}<x,\xi>}-1)\rho_1^{-H_0-d/2}(\xi)\mathrm{d}\widehat{W}(\xi),x\mapsto\int_{|\xi|\leq
R}(\mathrm{e}^{\mathrm{i}<x,\xi>}-1)\rho_2^{-H_0-d/2}(\xi)\mathrm{d}\widehat{W}(\xi)
\]
both belong to $\mathcal{C}^{r}(K,D_0)\hookrightarrow B^s_{p,q}(\mathbb{R}^d,D_0)$ for any compact subset $K$ of $\mathbb{R}^d$ and any
$r>\min_{\lambda\in Sp(\Delta)}\lambda$. Since any
$(\;E_0^{\!\!\!\!\!\!\!t},\mathbb{R}^d)$ pseudo--norms $\rho_1,\rho_2$ are
equivalent, Lemma~\ref{lem:eqPN} is then a straightforward
consequence of Theorem~1.1 of \cite{0} applied with $B=B^s_{p,q}(\mathbb{R}^d,D_0)$ which is either a separable Banach space either the dual of the separable space $B=B^{-s}_{p',q'}(\mathbb{R}^d,D_0)$ with $p',q'$ the respective conjugates of $p,q$ and to
\[
f_X=\rho_1^{-2H_0-d}\mathrm{1}_{|\xi|\geq R}\mbox{ and }\,f_Y=\rho_2^{-H_0-d}\mathrm{1}_{|\xi|\geq R}\;,
\]
and
\[
f_X=\rho_2^{-2H_0-d}\mathrm{1}_{|\xi|\geq R}\mbox{ and }\,f_Y=\rho_1^{-H_0-d}\mathrm{1}_{|\xi|\geq R}\;,
\]
successively.

Thus, using Lemmas~\ref{lem:cinfty} and~\ref{lem:eqPN}, we assume without loss of generality from now that the $(\mathbb{R}^d,E_0)$ pseudo--norm $\rho_{E_0}$, used to define the field $\{X_{\rho_{E_{0}},H_0}(x)\}_{x\in\mathbb{R}^{d}}$ belongs to $\mathcal{C}^{\infty}(\mathbb{R}^{d}\setminus\{0\})$. We shall use this assumption when proving that the wavelet coefficients of the field are weakly dependent (see Section~\ref{s:techlemmas}).

Observe that, to prove our local regularity results, we have to investigate the sample paths properties of $\varphi X$ for any function $\varphi\in\mathcal{D}(\mathbb{R}^d)$, that is for any $k_0\in\mathbb{Z}^d$, $r_0>0$ and any $\varphi\in\mathcal{D}(\mathbb{R}^d)$ satisfying $\mathrm{supp}(\varphi)\subset
B_{D_0}(k_0,r_0)$. Since the Besov spaces are invariant by translations
and dilatations, we may assume that $k_0=0$ and $r_0=1$. We have then to study the sample paths properties of the field $\varphi X$ for any function $\varphi\in\mathcal{D}(\mathbb{R}^d)$ such that $\mathrm{supp}(\varphi)\subset
B_{D_0}(0,1)$.

Our results come from the series expansion of
$X_{\rho_{E_{0}},H_0}$ in a Daubechies anisotropic wavelet basis
(see Section~\ref{SecBaseWAnis} just above). Recall that for any
$j\in \mathbb{N}, (G,\gamma)\in I_j(D)$, the wavelet coefficients
of $X_{\rho_{E_{0}},H_0}$ are defined as
\[
c_{j,G,\gamma}^{k}=<X_{\rho_{E_{0}},H_0},2^{\mathrm{Tr}(D_{j,G,\gamma})}\Psi_{j,G,\gamma}^{k}>\;.
\]
Fix now $|\cdot|_{D_0}$ a $(\mathbb{R}^{d},D_0)$ pseudo--norm.
Define $\Gamma_0(D_0)=\emptyset$ for $j=0$ and for any $j\geq 1$
\begin{equation}\label{e:defGammajD}
\Gamma_j(D_0)=\{k\in\mathbb{Z}^d,\,|k|_{D_0}< j2^{j}\}\;.
\end{equation}

Thereafter set
\begin{equation}\label{e:X1}
X_{\rho_{E_{0}},H_0}^{(1)}(x)=
\sum_{j,G,\gamma}\sum_{k\in\Gamma_j(D_0)}c_{j,G,\gamma}^{k}(\omega)\Psi_{j,G,\gamma}^{k}(x),
\end{equation}
and
\begin{equation}\label{e:X2}
X_{\rho_{E_{0}},H_0}^{(2)}(x)=\sum_{j,G,\gamma}\sum_{k\not\in\Gamma_j(D_0)}c_{j,G,\gamma}^{k}(\omega)\Psi_{j,G,\gamma}^{k}(x)\;.
\end{equation}
We will investigate separately the local sample path properties in anisotropic
Besov spaces of the two Gaussian fields $
X_{\rho_{E_{0}},H_0}^{(1)}$ and $
X_{\rho_{E_{0}},H_0}^{(2)}$. We first prove that
\begin{proposition}\label{pro:X1}Let $1\leq p,q\leq \infty$.
\begin{enumerate}
\item Almost surely, for any $\beta>1/q+d/\rho_{\min}(E_0)+\delta(p)$, the sample path of the field
$\{X_{\rho_{E_{0}},H_0}^{(1)}(x)\}_{x\in\mathbb{R}^{d}}$ belongs to $B^{H_0}_{p,q,|\log|^{\beta}}(\mathbb{R}^d,D_0)$.\\
\item Let $\varphi$ such that $\mathrm{supp}(\varphi)\subset B_{D_0}(0,1)$ and satisfying
\[
\varphi\equiv 1\mbox{ on }B_{D_0}(0,1/2)\;.
\]
Then almost surely, for $\beta=1/q+d/\rho_{\min}(E_0)+\delta(p)$ the sample path of the field
$\{\varphi X_{\rho_{E_{0}},H_0}^{(1)}(x)\}_{x\in\mathbb{R}^{d}}$ does not belong to $B^{H_0}_{p,q,|\log|^{-\beta}}(\mathbb{R}^d,D_0)$.
\end{enumerate}
\end{proposition}
{\bf Proof.} The proof uses several technics introduced
in~\cite{10}. The result comes from a comparison between
$\left[\sum_{k\in \Gamma_j(D_0)}|c_{j,G,\gamma}^k|^p\right]^{1/p}$
and $\left[\mathbb{E}(|c_{j,G,\gamma}^k|^2)\right]^{1/2}$
 and from Lemma~\ref{LemKey} which gives an estimate of $\left[\mathbb{E}(|c_{j,G,\gamma}^k|^2)\right]^{1/2}$. Set
\begin{equation}\label{eq:defgjk}
g_{j,G,\gamma}^{k}=\frac{c_{j,G,\gamma}^k}{\left[\mathbb{E}(|c_{j,G,\gamma}^k|^2)\right]^{1/2}}\;.
\end{equation}
for any $j\in\mathbb{N}$, $(G,\gamma)\in I_j$ and $k\in\Gamma_j(D_0)$. We need to distinguish two cases~: $p\neq \infty$ and $p=\infty$. In each case, we prove successively points~(i) and (ii).

Assume first that $p\neq \infty$ and let us prove point~(i) in this case. The definition of the sequence $(g_{j,G,\gamma}^{k})$ and the stationarity for any $(j,G,\gamma)$ of the sequence $(c_{j,G,\gamma}^{k},k\in\mathbb{Z}^d)$ implies that for any $j,G,\gamma$
\[
\left(\sum_{k\in\Gamma_j(D_0)}|c_{j,G,\gamma}^{k}|^p\right)^{1/p}
=\left[\mathbb{E}(|c_{j,G,\gamma}^0|^2)\right]^{1/2}\cdot
\left(\sum_{k\in\Gamma_j(D_0)}|g_{j,G,\gamma}^{k}|^p\right)^{1/p}\;.
\]
Use now the weak correlation of the wavelet coefficients and the
two estimates of
$\left[\mathbb{E}(|c_{j,G,\gamma}^0|^2)\right]^{1/2}$ and
of $n_j=\mathrm{card}(\Gamma_j(D_0))$ respectively proved in
Lemmas~\ref{LemKey} and~\ref{lem:estimnj}. One deduces that the
following inequality holds for any $j\geq 0$~:
\begin{equation}\label{e:ineq1}
\left(\sum_{k\in\Gamma_j(D_0)}|c_{j,G,\gamma}^{k}|^p\right)^{\frac{1}{p}}\leq
C
2^{j(\frac{d}{p}-H_0)}j^{d^*}\left(\frac{1}{n_{j}}\sum_{k\in\Gamma_j(D_0)}|g_{j,G,\gamma}^{k}|^p\right)^{\frac{1}{p}}\;,
\end{equation}
where
\[
d^*=\frac{d}{2\rho_{\min}(E_0)}+\frac{d}{p}\;.
\]
Lemma~\ref{lem:CRK1} stating a central limit theorem for the sequence $(g_{j,G,\gamma}^{k})$
and inequality~(\ref{e:ineq1}) then prove point~(i) of the proposition for the case $p<\infty$.

We now prove point~(ii) for $p\neq \infty$. Set now
\[
\Gamma_j'(D_0)=\{k\in\mathbb{Z}^d,\,|k|_{D_0}\leq 2^j/j\}\;.
\]
Using the assumptions on the support of $\varphi$, remark that for $j$ sufficiently large and for any $k\in\Gamma_j'(D_0)$, one has
\[
c_{j,G,\gamma}^k(\varphi X)=c_{j,G,\gamma}^k(X)\;.
\]
Use the same arguments as in the proof of Lemma~\ref{lem:estimnj} and deduce that $n'_{j}=\mathrm{Card}(\Gamma'_j(D_0))\sim j^{-d}2^{jd}$. Since $\Gamma'_j(D_0)\subset \Gamma_j(D_0)$, a similar approach to above then yields that for some $C>0$ and for any $j\geq 1$
\[
\left(\sum_{k\in\Gamma_j(D_0)}|c_{j,G,\gamma}^{k}|^p\right)^{\frac{1}{p}}\geq \left(\sum_{k\in\Gamma_j'(D_0)}|c_{j,G,\gamma}^{k}|^p\right)^{\frac{1}{p}}\geq
C 2^{j(\frac{d}{p}-H_0)}j^{-d^*}\left(\frac{1}{n'_{j}}\sum_{k\in\Gamma_j'(D_0)}|g_{j,G,\gamma}^{k}|^p\right)^{\frac{1}{p}}\;.
\]
which directly implies point~(ii) of the proposition.

If $p=\infty$, a similar approach implies that almost surely there
exists some $C_1,C_2>0$ such that for any $j,G,\gamma$
\[
C_1 2^{-jH_0}j^{-d^*}\left(\frac{1}{\sqrt{\log(n_{j})}}\sup_{k\in\Gamma_j}|g_{j,G,\gamma}^{k}|\right)\leq\left(\sum_{k\in\Gamma_j}|c_{j,G,\gamma}^{k}|^p\right)^{\frac{1}{p}}\;,
\]
and
\[
\left(\sum_{k\in\Gamma_j}|c_{j,G,\gamma}^{k}|^p\right)^{\frac{1}{p}}\leq C_2 2^{-jH_0}j^{d^*}\left(\frac{1}{\sqrt{\log(n_{j})}}\sup_{k\in\Gamma_j}|g_{j,G,\gamma}^{k}|\right)\;,
\]
with
\[
d^*=\frac{d}{2\rho_{\min}(E_0)}\;.
\]
Lemma~\ref{lem:CRK2} and the inequality just above then implies the result stated in point~(i) for the case $p=\infty$. The proof of point~(ii) for $p=\infty$ also follows from the above inequality replacing $\Gamma_j(D_0)$ with $\Gamma'_j(D_0)$ as in the case $p\neq \infty$.

We now investigate the sample paths properties of
$\{X_{\rho_{E_{0}},H_0}^{(2)}(x)\}_{x\in\mathbb{R}^{d}}$.
\begin{proposition}\label{pro:X2}
Almost surely, the sample path of the field
$\{X_{\rho_{E_{0}},H_0}^{(2)}(x)\}_{x\in\mathbb{R}^{d}}$
belong to  $B^{H'}_{p,q, loc}(\mathbb{R}^{d},E_0)$ for any
\[
0<H_0<H'<\rho_{\min}(D_0)=\rho_{\min}(E_0)\;,
\]
and any $1\leq p,q\leq \infty$.
\end{proposition}
{\bf Proof.} Using the transference results of~\cite{36} (see Theorem 5.28) and the usual embedding of isotropic Besov spaces defined on bounded domains one remarks that
\[
\mathcal{C}^{s+\varepsilon}_{loc}(\mathbb{R}^{d},D_0)\subset
B^{s}_{p,q, loc}(\mathbb{R}^{d},D_0)\;,
\]
for any $1\leq p,q\leq \infty$ and any $s,\varepsilon>0$. It then suffices to prove the result for $p=q=\infty$.\\

Let now consider $H'\in (H_0,\rho_{\min}(E_0))$, $\varepsilon>0$
and $\varphi\in\mathcal{D}(\mathbb{R}^{d})$. Recall that we
assumed that
\[
\mathrm{supp}(\varphi)\subset B_{D_0}(0,1)=\{x,|x|_{D_0}\leq
1\}\;,
\]
and $0\leq\varphi\leq 1$ on $\mathbb{R}^{d}$. We denote by $Y$ the random field $\varphi
X_{\rho_{E_{0}},H_0}^{(2)}$.\\

We will give an upper bound of $|Y(x+h)-Y(x)|$ for any given $x$
in $B_{D_0}(0,1)$ and $h$ sufficiently small. Observe that
\[
Y(x+h)-Y(x)=Y_1(x,h)+Y_2(x,h)\;,
\]
with
\begin{eqnarray*}
Y_1(x,h)=\sum\limits_{j,G,\gamma}\sum_{k\not\in\Gamma_j(D_0)}c_{j,G,\gamma}^{k}(\varphi(x+h)-\varphi(x))\Psi_{j,G,\gamma}^{k}(x)\;,\\
Y_2(x,h)=\sum\limits_{j,G,\gamma}\sum_{k\not\in\Gamma_j(D_0)}c_{j,G,\gamma}^{k}\varphi(x+h)(\Psi_{j,G,\gamma}^{k}(x+h)-\Psi_{j,G,\gamma}^{k}(x))\;.
\end{eqnarray*}

We first bound $Y_1(x,h)$. Let
$\varepsilon=1-H'/\rho_{\min}(E_0)$. We now use that $\varphi\in
B_{\infty,\infty,loc}^{1-\varepsilon}(\mathbb{R}^d)$, $h$
sufficiently small and $x$ belongs to the compact set
$B_{D_0}(0,1)$. Hence , by Lemma~\ref{lem:CRK2} and the fast decay
of the wavelets, almost surely for any $M>0$ and for some $C>0$
one has
\[
|Y_1(x,h)| \leq C|h|^{1-\varepsilon}
\sum_{j,G,\gamma}j^{d^*}2^{-jH_0}
\left(\sum\limits_{k\not\in\Gamma_j(D_0)}\frac{1}{(1+|k-2^{D_{j,G,\gamma}}x|)^M}\right)\;.
\]
Here we denoted $d^*=1/2+d/\rho_{\min}(E_0)$. Further, by
assumption on $k$ and $x$
\[
|k|_{D_0}\geq j2^j\geq j|2^{D_{j,G,\gamma}}x|_{D_0}\;.
\]
Since $|\cdot|_{D_0}$ is a $(\mathbb{R}^{d},D_0)$ pseudo--norm, by
the triangular inequality~(\ref{e:triangineq}), one deduces that
for $j$ sufficiently large
\[
|k-2^{D_{j,G,\gamma}}x|_{D_0}\geq C|k|_{D_0}\;.
\]
for some $C\in (0,1)$. Then by comparison between $|\cdot|_{D_0}$
and the usual Eucidean norm, one deduces that there exists some
$\alpha>0$ such that for $j$ sufficiently large and any $x$ in
$B_{D_0}(0,1)$
\[
|k-2^{D_{j,G,\gamma}}x|\geq (|k|/2)^{\alpha}\;.
\]
Then
\[
|Y_1(x,h)| \leq
C|h|^{1-\varepsilon}\left(\sum\limits_{j,G,\gamma}j^{d^*}2^{-jH_0}\sum_{k\not\in
\Gamma_j(D_0)}\frac{1}{(1+|k|^{\alpha})^M}\right)\;.
\]
Since, for $M$ sufficiently large
\[
\sum\limits_{j,G,\gamma}j^{d^*}2^{-jH_0}\sum_{k\not\in
\Gamma_j(D_0)}\frac{1}{(1+|k|^{\alpha})^M}<\infty\;,
\]
one has almost surely $|Y_1(x,h)| \leq C'|h|_{D_0}^{H'}$. \\

By the same approach, we can bound $Y_2(x,h)$. Indeed, using the
fact that $\varphi$ is bounded and the mean value theorem for
$\Psi_{j,G,\gamma}^k$, we then prove that almost surely for some
$C>0$
\[
|Y_2(x,h)|\leq
C\sum_{j,G,\gamma}j^{d^*}2^{-jH_0}|2^{D_{j,G,\gamma}}h|
\left(\sup_{y\in
[x,x+h]}\sum_{k\not\in\Gamma_j(D_0)}\frac{1}{(1+|k-2^{D_{j,G,\gamma}}y|)^M}\right)\;.
\]
The end of the proof is exactly the same as above remarking that
\[
|2^{D_{j,G,\gamma}}h|\leq
j^{\delta}2^j|h|_{D_0}^{\rho_{\min}(E_0)}\;,
\]
for some $\delta>0$. Proposition~\ref{PropLocRegAnisEgale} then follows directly from
Propositions~\ref{pro:X1} and~\ref{pro:X2}.
\subsection{Proof of regularity results in anisotropic Besov spaces with an anisotropy commuting with this of the field}
The following proposition extends the results of
Proposition~\ref{PropLocRegAnisEgale} in anisotropic Besov spaces
$B^s_{p,q}(\mathbb{R}^d,E)$ with $E\in\mathcal{E}^+_d$ commuting
with $E_0$.
\begin{proposition}\label{PropRegBesovAutreAnis}Let $1\leq p,q \leq +\infty$, $\varepsilon>0$ and $E\in\mathcal{E}^+_d$
commuting with $E_0$. Then
\begin{enumerate}
\item Almost surely the sample path of
$\{X_{\rho_{E_0},H_0}(x)\}_{x\in\mathbb{R}^d}$ belongs to
$B^{H_0\frac{\lambda_{m}}{\lambda_{m}^{0}}-\varepsilon}_{p,q,\mathrm{loc}}(\mathbb{R}^{d},E)$.
\item Almost surely the sample path of
$\{X_{\rho_{E_0},H_0}(x)\}_{x\in\mathbb{R}^d}$ does not belong to
$B^{H_0\frac{\lambda_{m}}{\lambda_{m}^{0}}+\varepsilon}_{p,q,\mathrm{loc}}(\mathbb{R}^{d},E)$.
\end{enumerate}
\end{proposition}
The proof is made in several steps. First we need to compare Besov
spaces with different commuting anisotropies.
\subsubsection{A comparison result between Besov spaces with different commuting anisotropies}
Since $E$ and $E_0$ are commuting, we can then assume (up to a
change of basis) that $D_0$ and $D$ are two diagonal matrices of
the form :
\begin{equation}\label{EqMatricesDDp}
D_0=\begin{pmatrix}\lambda_{1}^{0}Id_{d_{1}}&
&0\\&\ddots&\\0&&\lambda_{m}^{0}Id_{d_{m}}\end{pmatrix},\,D=\begin{pmatrix}\lambda_{1}Id_{d_{1}}&
&0\\&\ddots&\\0&&\lambda_{m}Id_{d_{m}}\end{pmatrix},
\end{equation}
with
\begin{equation}\label{EqValeursPropres}
\frac{\lambda_{m}}{\lambda_{m}^{0}}\leq\cdots\leq
\frac{\lambda_{1}}{\lambda_{1}^{0}}\;.
\end{equation}
\begin{proposition}\label{PropCompBesAnis2}
The notations and assumptions are as above.
For any $\alpha>0$, $\beta\in\mathbb{R}$ and $p,q\in (0,+\infty]$, one has the following embedding
\[
B^{\alpha}_{p,q,|\log|^{\beta}}(\mathbb{R}^d,D_0) \hookrightarrow
B^{\alpha\frac{\lambda_{m}}{\lambda_{m}^{0}}}_{p,q,|\log|^{\beta}}(\mathbb{R}^d,D)\;.
\]
\end{proposition}
The proof is straightforward and based on finite differences
characterization of Besov spaces given in Theorem 5.8~(ii)
of~\cite{36}.
\subsubsection{Proof of Proposition~\ref{PropRegBesovAutreAnis}}
We only prove the second point of
Proposition~\ref{PropRegBesovAutreAnis} since the first one is a
straigthforward consequence of Propositions~\ref{PropCompBesAnis1}
and~\ref{PropCompBesAnis2}.  To this end we use the following
characterization of anisotropic Besov spaces
$B^s(\mathbb{R}^{d},\Delta)$ with diagonal anisotropy $\Delta$ (see Theorem~5.8 of \cite{36})~:
\begin{proposition}\label{pro:characBesov}
Let $\Delta$ a matrix belonging to $\mathcal{E}^+$ of the form
\[
\Delta=\begin{pmatrix}\alpha_1&0&\cdots&0\\0&\alpha_2&\hdots&0\\\vdots&\vdots&\ddots&\vdots\\0&0&\hdots&\alpha_d\end{pmatrix}\;.
\]
and $s\in (0,\rho_{\min}(\Delta))$, $M_\ell=[s/\alpha_\ell]+1$ for any $\ell=1,\cdots,d$. Then $f\in
B^{s}_{p,p}(\mathbb{R}^{d},\Delta)$ if and only if
\[
\|f\|_{L^p}+\sum_{\ell=1}^d\left(\int_{0}^1 \|(\Delta^{M_\ell}_{t e_\ell} f)(x)\|_{L^p}^p
t^{-sp/\alpha_{\ell}-1}\rmd t\right)^{1/p}<\infty\;,
\]
where $(e_\ell)$ is the canonical basis of $\mathbb{R}^{d}$ and where as usual, for any $x,h\in\mathbb{R}^d$
\[
(\Delta_h^1 f)(x)=f(x+h)-f(x),\,\cdots,\,(\Delta_h^{M_\ell} f)(x)=(\Delta_h^1 \Delta^{M_{\ell}-1}f)(x)\;.
\]
\end{proposition}
\begin{remark}
If for $\ell=1,\cdots,d$,$s\in (0,\alpha_\ell)$ then $M_\ell=1$.
\end{remark}
{\bf Proof.} This statement is proved in Theorem~5.8 in~\cite{36}.

We now prove Proposition~\ref{PropRegBesovAutreAnis}. We first
remark that we have only to consider the case where $E_0=D_0$.
Indeed, let $|\cdot|_{E_0}$ (resp $|\cdot|_{D_0}$) be a
$(\mathbb{R}^{d},E_0)$ (resp a $(\mathbb{R}^{d},D_0)$)
pseudo--norm. Lemma~\ref{LemComPNDiagCom} then implies that for
any $\varepsilon>0$ and any $|\xi|$ sufficiently large
\[
|\xi|_{D_0}^{1-\varepsilon}\leq |\xi|_{E_0}\leq
|\xi|_{D_0}^{1+\varepsilon}\;.
\]
Hence Theorem~1.1 of~\cite{0} applied successively with $f_X=\mathrm{1}_{|\xi|\leq R}|\xi|_{E_0}^{-2H_0-d}, f_Y=\mathrm{1}_{|\xi|\leq R}|\xi|_{D_0}^{-(1-\varepsilon)(2H_0+d)}$,
$f_X=\mathrm{1}_{|\xi|\leq R}|\xi|_{D_0}^{-(1+\varepsilon)(2H_0+d)},f_Y=\mathrm{1}_{|\xi|\leq R}|\xi|_{E_0}^{-2H_0-d}$ as above and Proposition~\ref{PropRegBesovAutreAnis} proved in the case $E_0=D_0$ yield the result in the general case $E_0\in\mathcal{E}^+_d$.

From now, we then assume that $E_0=D_0$ and that the
$(\mathbb{R}^{d},D_0)$ pseudo--norm involved in the construction
of the studied field belongs to
$\mathcal{C}^\infty(\mathbb{R}^{d}\setminus\{0\})$ which ensures the weak correlation of the wavelet coefficients.

As in the proof of Proposition~\ref{PropLocRegAnisEgale}, we use
an expansion of the Gaussian field $X_{\rho_{D_0},H_0}$ in a
Daubechies wavelet basis and we define,
\begin{equation}\label{e:X1t}
\widetilde{X}^{(1)}(x)=
\sum_{j,G,\gamma}\sum_{k\in\widetilde{\Gamma}_j(D_0)}c_{j,G,\gamma}^{k}(\omega)\Psi_{j,G,\gamma}^{k}(x),
\end{equation}
and
\begin{equation}\label{e:X2t}
\widetilde{X}^{(2)}(x)=\sum_{j,G,\gamma}\sum_{k\not\in\widetilde{\Gamma}_j(D_0)}c_{j,G,\gamma}^{k}(\omega)\Psi_{j,G,\gamma}^{k}(x)\;.
\end{equation}
where $\widetilde{\Gamma}_j(D_0)=\{k\in\mathbb{Z}^d,\, |k|_{D_0}\leq
C_1 2^j\}$.

As in the proof of Proposition~\ref{PropLocRegAnisEgale}, we see
that, for $C_1$ sufficiently large, almost surely
$\widetilde{X}^{(2)}$ belongs to
$B^{H_0\lambda_m/\lambda_m^0+\varepsilon}_{p,q,loc}(\mathbb{R}^{d},D)$
for any $1\leq p,q\leq \infty$ and $\varepsilon>0$ sufficiently small.

We then have to prove our a.s. non local regularity results for the Gaussian field
$\widetilde{X}^{(1)}$. Remark now that since the multiresolution analysis is compactly supported so is $\widetilde{X}^{(1)}$.  To show point~(ii) of Proposition~\ref{PropRegBesovAutreAnis}, it is then sufficient to prove that a.s. the sample paths of $\widetilde{X}^{(1)}$ does not belong to $B^{\lambda_m H_0/\lambda^0_m+\varepsilon}(\mathbb{R}^d,D)$ for any $\varepsilon>0$.

Set $M=[s/\lambda_m]+1$ which may be greater than one. In view of Proposition~\ref{pro:characBesov}, we shall then give an almost sure lower bound of
\[
\|(\Delta^M_{te_d} \widetilde{X}^{(1)})(x)\|_{L^p}^p=\int_{\mathbb{R}^{d}}|
(\Delta^M_{te_d} \widetilde{X}^{(1)})(x)|^p\mathrm{d}x\;,
\]
for any $p\geq 1$ and any $t$ of the form
$t=2^{-[j_0\lambda^0_m]}$ where $j_0$ is a fixed non--negative integer.

Set
$$\Delta_{j_0}=\begin{pmatrix}[j_0\lambda^0_1]&&0\\
&\ddots&\\
0&&[j_0\lambda^0_m]\end{pmatrix}\;.$$
Observe that if $t=2^{-[j_0\lambda^0_m]}$, one has $t e_d=2^{-\Delta_{j_0}} e_d$. Remark also that $\widetilde{X}^{(1)}$ can be written as the sum of its low frequency component and its high frequency component, namely that
\[
\widetilde{X}^{(1)}=
\widetilde{X}^{(1)}_{LF}+\widetilde{X}^{(1)}_{HF}
\]
with
\[
\widetilde{X}^{(1)}_{LF}(x)=\sum_{j\leq j_0}\sum_{G,\gamma}\sum_{k\in\widetilde{\Gamma}_j(D_0)
}c_{j,G,\gamma}^k\Psi_{j,G,\gamma}^k(x)\mbox{ and }
\widetilde{X}^{(1)}_{HF}(x)=\sum_{j\geq j_0+1
}\sum_{G,\gamma}\sum_{k\in\widetilde{\Gamma}_j(D_0)
}c_{j,G,\gamma}^k\Psi_{j,G,\gamma}^k(x)\;.
\]
Using the triangular inequality, one has
\begin{equation}\label{e:triangineqLp}
\|
(\Delta^M_{te_d} \widetilde{X}^{(1)})(x)\|_{L^p}\geq \|
(\Delta^M_{te_d} \widetilde{X}^{(1)}_{HF})(x)\|_{L^p}-\|
(\Delta^M_{te_d} \widetilde{X}^{(1)}_{LF})(x)\|_{L^p}
\end{equation}
To give a lower bound of $\|
(\Delta^M_{te_d} \widetilde{X}^{(1)})(x)\|_{L^p}^p$, we shall then give a lower bound of $\|
(\Delta^M_{te_d} \widetilde{X}^{(1)}_{HF})(x)\|_{L^p}^p$ and an upper bound of $\|
(\Delta^M_{te_d} \widetilde{X}^{(1)}_{LF})(x)\|_{L^p}^p$.

Let us first give an upper bound of $\|
(\Delta^M_{te_d} \widetilde{X}^{(1)}_{LF})(x)\|_{L^p}^p$. We suppose that the multiresolution analysis is $s$ smooth for some $s\in (H_0,\rho_{\min}(D_0))$. By the finite differences definition of the spaces $\dot{B}^s_{p,\infty}(\mathbb{R}^d,D_0)$ and the fact that for any $M\geq 1$ $|\Delta^M_h f(x)|\leq \sum_{\ell=1}^M |f(x+\ell h)-f(x+(\ell-1)h)|$, one has
\begin{eqnarray*}
&&\left|\left|\sum_{j\leq j_0
}\sum_{G,\gamma}\sum_{k\in\widetilde{\Gamma}_j(D_0)
}c_{j,G,\gamma}^k\left(\Delta^M_{2^{-\Delta_{j_0}}e_d}\Psi_{j,G,\gamma}^k\right)(x)\right|\right|_{L^p}\\
&\leq&C\sum_{j\leq j_0
}|2^{-\Delta_{j_0}}e_d|_{D_0}^s
\left|\left|\sum_{G,\gamma}\sum_{k\in\widetilde{\Gamma}_j(D_0)}c_{j,G,\gamma}^k\Psi_{j,G,\gamma}^k\right|\right|_{\dot{B}^s_{p,\infty}(\mathbb{R}^d,D_0)}
\end{eqnarray*}
Use now the wavelet characterization of the homogeneous Besov spaces $\dot{B}^s_{p,\infty}(\mathbb{R}^d,D_0)$. Then for some $C>0$
\[
\left|\left|\sum_{G,\gamma}\sum_{k\in\widetilde{\Gamma}_j(D_0)}c_{j,G,\gamma}^k\Psi_{j,G,\gamma}^k\right|\right|_{\dot{B}^s_{p,\infty}(\mathbb{R}^d,D_0)}\leq C2^{j(s-d/p)}
\left(\sum_{G,\gamma}\sum_{k\in\widetilde{\Gamma}_j(D_0)}|c_{j,G,\gamma}^k|^p\right)^{1/p}\;.
\]
As the proof of Proposition~\ref{PropLocRegAnisEgale}, we can estimate a.s. $\sum_{k\in\widetilde{\Gamma}_j(D_0)}|c_{j,G,\gamma}^k|^p$. Hence, we
deduce that there exists an a.s. positive constant $C'$ such that
\begin{equation}\label{e:boundLF}
\left|\left|
\left(\Delta^M_{t e_{d}}\widetilde{X}^{(1)}_{LF}\right)(x)\right|\right|_{L^p(\mathbb{R}^{d})}^p\leq C2^{-j_0 s}\sum_{j\leq j_0}2^{j(s-\frac{d}{p})}\cdot\left(2^{j\frac{d}{p}}j^{\frac{d}{2\rho_{\min}(E_0)}}2^{-jH_0}\right)\leq C'\;,
\end{equation}
where $C'$ is not depending on $j_0$ nor $s$.

We now give a lower bound of $\|
(\Delta^M_{t e_{d}}\widetilde{X}^{(1)}_{HF})(x)\|_{L^p}^p$. To this end perform the change of variable $x=2 ^{-\Delta_{j_0}}y$ and deduce that
\begin{equation}\label{e:chgvar}
\|
(\Delta^M_{t e_{d}}\widetilde{X}^{(1)}_{HF})(x)\|_{L^p}^p=2^{-\mathrm{Tr}(\Delta_{j_0})}\|
(\Delta^M_{t e_{d}}\widetilde{X}^{(1)}_{HF})(2 ^{-\Delta_{j_0}}y)\|_{L^p}^p
\end{equation}
By definition of $\widetilde{X}^{(1)}_{HF}$ one has
\[
(\Delta^M_{t e_{d}}\widetilde{X}^{(1)}_{HF})(2 ^{-\Delta_{j_0}}y)=\sum_{j\geq j_0+1
}\sum_{G,\gamma}\sum_{k\in\widetilde{\Gamma}_j(D_0)
}c_{j,G,\gamma}^k\left(\Delta^M_{t e_{d}}\Psi_{j,G,\gamma}^k\right)(2 ^{-\Delta_{j_0}}y)\;.
\]
Define now for any $j'\geq 1$, any $G\in
\{F,M\}^{d^*}$ and any $\gamma'\in \mathbb{N}^d$ such that
\[
[(j'-1)\lambda_r]-2\leq \gamma'_r\leq [j'\lambda_r]+2
\]
the family of functions
\[
h_{j',G,\gamma'}^k(y)=(\Delta^M_{e_d}\Psi^{(G)})(2^{D_{j',G,\gamma'}}y-k)\;,
\]
where
\[
D_{j',G,\gamma'}=\begin{pmatrix}\gamma'_1&0&\hdots&0\\0&\ddots&&\vdots\\\vdots&\ddots&\gamma'_{m-1}&0\\0&\hdots&0&\gamma'_m\end{pmatrix}\;,
\]
and for $j'= 0$ and any $k\in\mathbb{Z}^d$
$h_{j',(F,\cdots,F),(0,\cdots,0)}^k(x)=\Psi^{(F,\cdots,F)}(x-k)$. Observe that this is a family of
inhomogeneous smooth analysis molecules in the sense of
Definition~5.3 in~\cite{8}. Further, if $j'=j-j_0$ and $\widetilde{\gamma}=([j_0 \lambda_r^0])_{r}$ one has
\begin{eqnarray*}
\widetilde{X}^{(1)}_{HF}(2^{-\Delta_{j_0}}(y+
e_{d}))-\widetilde{X}^{(1)}_{HF}(2^{-\Delta_{j_0}}y)=\sum_{j'\geq 1}\sum_{G,\gamma'}\sum_{k\in\widetilde{\Gamma}_j(D_0)
}\widetilde{c}_{j',G,\gamma'}^k h_{j',G,\gamma'}^k(y)\;,
\end{eqnarray*}
with $\widetilde{c}_{j',G,\gamma'}^k=c_{j'+j_0,G,\gamma'+\widetilde{\gamma}}^k$ if $(G,\gamma'+\widetilde{\gamma})\in I^{j'+j_0}(D_0)$ and $\widetilde{c}_{j',G,\gamma}^k=0$ otherwise. Hence
\begin{eqnarray*}
\|(\Delta^M_{2^{-\Delta_{j_0}} e_d}\widetilde{X}^{(1)}_{HF})(2^{-\Delta_{j_0}}y)\|_{L^{p}}&\geq&\|(\Delta^M_{2^{-\Delta_{j_0}} e_d}\widetilde{X}^{(1)}_{HF})(2^{-\Delta_{j_0}}y)\|_{B^0_{p,p,|\log|}(\mathbb{R}^{d},D_0)}\\
&\geq&\left(\sum_{j'\geq 1}\sum_{(G,\gamma)\in I^{j'+j_0}}j' 2^{-j'd}\sum_{k\in
\widetilde{\Gamma}_{j'+j_0}(D_0)}|c_{j'+j_0,G,\gamma}^k|^p\right)^{1/p}
\end{eqnarray*}
We use once more an a.s. estimate of $\sum_k|c_{j,G,\gamma}^k|^p$ as
in the proof of Proposition~\ref{PropLocRegAnisEgale}. Since as in Lemma~\ref{lem:estimnj}, we can prove that $\mathrm{Card}(\widetilde{\Gamma}_{j'+j_0}(D_0)\geq 2^{(j'+j_0)d}$. Hence there exists an a.s. positive constant $C$ such that
\[
\sum_{k\in
\widetilde{\Gamma}_{j'+j_0}(D_0)}|c_{j'+j_0,G,\gamma}^k|^p\geq 2^{-(j'+j_0)(H_0 p-d)}\;.
\]
Hence
\[
\|(\Delta^M_{2^{-\Delta_{j_0}} e_d}\widetilde{X}^{(1)}_{HF})(2^{-\Delta_{j_0}}y)\|_{L^p}^p\geq
C\sum_{j',G,\gamma}j'
2^{-j'd}2^{-(j'+j_0)(H_0 p-d)}
\]
Use now the last inequality and relation~(\ref{e:chgvar}). Then a.s.
\[
\|(\Delta^M_{2^{-\Delta_{j_0}} e_d}\widetilde{X}^{(1)}_{HF})(2^{-\Delta_{j_0}}y)\|_{L^p}^p\geq
C2^{-j_0 d}j_0^{-d}\left(\sum_{j'\geq 0}j'
2^{-j'd}2^{-(j'+j_0)(H_0 p-d)}\right)\geq C2^{-j_0 H_0 p}\;.
\]
We deduce that a.s.
\begin{eqnarray*}
&&\int_0^1 \| (\Delta^M_{2^{-[j_0\lambda^0_m]}
e_{\ell}}\widetilde{X}^{(1)}_{HF})(x)\|_{L^p}^p t^{-s/\lambda_m
p-1}\mathrm{d}t\\
&\geq&\sum_{j_0=0}^{+\infty}\int_{2^{-[(j_0+1)\lambda^0_m]}}^{2^{-[j_0\lambda^0_m]}} \|(\Delta^M_{2^{-[j_0\lambda^0_m]}
e_{\ell}}\varphi \widetilde{X}^{(1)}_{\rho_{D_0},H_0})(x)\|_{L^p(\mathbb{R}^{m})}^p t^{-s/\lambda_m
p-1}\mathrm{d}t\\
&\geq& C\sum_{j_0=0}^{+\infty}j_0^{-d} 2^{-j_0 H_0
p}\left(2^{-[j_0\lambda^0_m]}\right)^{-sp/\lambda_m-1}2^{-[j_0\lambda^0_m]}\\
&=&C\sum_{j_0=0}^{+\infty}2^{-j_0 H_0 p}j_0^{-d}
\left(2^{-[j_0\lambda^0_m]}\right)^{-sp/\lambda_m}\;.
\end{eqnarray*}
Since
\[
\sum_{j_0=0}^{+\infty}j_0^d 2^{-j_0 H_0
p}\left(2^{-[j_0\lambda^0_m]}\right)^{-sp/\lambda_m}=+\infty\;,
\]
if
$s\lambda^0_m/\lambda_m-H_0 p> 0$, we deduce that a.s.
\[
\int_0^1 \| \widetilde{X}^{(1)}_{HF}(x+2^{-[j_0\lambda^0_m]}
e_{\ell})-\widetilde{X}^{(1)}_{HF}(x)\|_{L^p}^pt^{-s/\lambda_m
p-1}\mathrm{d}t=+\infty
\]
for $s> H_0\lambda_m/\lambda^0_m$. Using (\ref{e:boundLF}) and the triangular inequality~(\ref{e:triangineqLp}), it ends the proof of
Proposition~\ref{PropRegBesovAutreAnis}. It also implies directly Theorem~\ref{ThOptim2}
\section{Technical lemmas}\label{s:techlemmas}
Our results about smoothness of the sample path are based on the following lemma
\begin{lemma}\label{LemKey}Assume that the anisotropic
multi--resolution analysis considered is $\mathcal{C}^1$ and admits at least one
vanishing moment.

Let $\{X_{\rho_{E_0},H_0}(x)\}_{x\in\mathbb{R}^{d}}$ the Gaussian field defined by~(\ref{EqDefX}) with $\rho=\rho_{E_0}$. Assume also that the pseudo--norm $\rho_{E_0}$ involved in the construction of this field is at least $\mathcal{C}^1(\mathbb{R}^d\setminus\{0\})$. Then the wavelet coefficients of the random
field $\{X_{\rho_{E_0},H_0}(x)\}_{x\in\mathbb{R}^{d}}$ are weakly
dependent in the following sense
\begin{enumerate}
\item There exists some $C_0>0$ such  for any $j\geq 1$, $(G,\gamma)\in I_j$ and $(k,k')\in (\mathbb{Z}^{d})^{2}$
\begin{equation}\label{eq:WC1}
|\mathbb{E}(c_{j,G,\gamma}^{k}c_{j,G,\gamma}^{k'})|\leq C_{0}\frac{j^{2d/\rho_{\min}} 2^{-2jH_0}}{1+|k-k'|}.
\end{equation}
\item There exists some $C_1,C_2>0$ such that for any $j\geq 1$, $(G,\gamma)\in I_j$ and any $k\in (\mathbb{Z}^{d})$
\begin{equation}\label{eq:WC2}
C_{1}j^{-d/\rho_{\min}(E_0)}2^{-2jH_{0}}\leq
\mathbb{E}(|c_{j,G,\gamma}^{k}|^{2})\leq
C_{2}j^{d/\rho_{\min}(E_0)}2^{-2jH_{0}}\;.
\end{equation}
\end{enumerate}
\end{lemma}
\begin{remark}
Theorem~1.1 of~\cite{CV13b} imply that, studying the sample paths properties of OSSRGF, we can always assume that the pseudo--norm $\rho_{E_0}$ involved in the construction of this field is $\mathcal{C}^\infty(\mathbb{R}^d\setminus\{0\})$. Then the assumptions of Lemma~\ref{LemKey} are satisfied.
\end{remark}
{\bf Proof.}Since the anisotropic multiresolution analysis admits at least one vanishing moment, one has $\widehat{\psi}_M(0)=0$. Further, for any $j\geq 1$, $(G,\gamma)\in I_j$ and for all $k\in \mathbb{Z}^{d}$
\[
c_{j,G,\gamma}^{k}=\int_{\mathbb{R}^{d}}\rme^{\rmi 2^{-\;D^{\!\!\!\!\!\!t}_{j,G,\gamma}}k\;\xi}
\overline{\widehat{\psi}^{(G)}}(2^{-\;D^{\!\!\!\!\!\!\!t}_{j,G,\gamma}}\;\xi)\rho_{\;E_0^{\!\!\!\!\!\!\!t}}\;(\xi)^{-H_0-d/2}\rmd\widehat{W}(\xi)\;.
\]
This formula implies that (set $\zeta=2^{-\;D^{\!\!\!\!\!\!\!t}_{j,G,\gamma}}\;\xi$)
\[
\mathbb{E}(|c_{j,G,\gamma}^{k}|^2)= 2^{j\mathrm{Tr}(D_{j,G,\gamma})}\int_{\mathbb{R}^{d}}
|\widehat{\psi}^{(G)}(\zeta)|^2\rho_{\;E_0^{\!\!\!\!\!\!\!\!t}}(2^{\;D^{\!\!\!\!\!\!\!t}_{j,G,\gamma}}\;\zeta)^{-2H_0-d}\rmd\zeta\;.
\]
Since $2^{(j-2)d}\leq \mathrm{Tr}(D_{j,G,\gamma})\leq 2^{j d}$, using Lemma~\ref{LemComPNDiagCom} and the inequalities $C_1 2^j\leq |2^{\;D^{\!\!\!\!\!\!\!t}_{j,G,\gamma}}\;\zeta|_{D_0}\leq C_2 2^j$ imply that
\[
\mathbb{E}(|c_{j,G,\gamma}^{k}|^2)\geq C_1 2^{-2j(H_0+d)}2^{jd}\int_{\mathbb{R}^{d}}
|\widehat{\psi}^{(G)}(\zeta)|^2|\zeta|_{\;D_0^{\!\!\!\!\!\!\!t}}^{-2H_0-d}\;(1+\log(|\zeta|_{\;D_0^{\!\!\!\!\!\!t}})+j)^{-d/\rho_{\min}(E_0)}\;\rmd\zeta\;,
\]
and
\[
\mathbb{E}(|c_{j,G,\gamma}^{k}|^2)\leq C_2 2^{-2j(H_0+d)}2^{jd}\int_{\mathbb{R}^{d}}
|\widehat{\psi}^{(G)}(\zeta)|^2|\zeta|_{\;D_0^{\!\!\!\!\!\!\!\!t}}^{-2H_0-d}\;(1+\log(|\zeta|_{\;D_0^{\!\!\!\!\!\!\!\!t}}\;)+j)^{d/\rho_{\min}(E_0)}\rmd\zeta\;.
\]
We then proved inequalities~(\ref{eq:WC2}).\\
To prove inequality~(\ref{eq:WC1}) remark that for any $\ell\in\{1,\cdots,d\}$
\[
(k_{\ell}-k'_{\ell})\mathbb{E}(c_{j,G,\gamma}^{k}c_{j,G,\gamma}^{k'})=\int_{\mathbb{R}^{d}}(k_{\ell}-k'_{\ell})\rme^{\rmi 2^{-\;D^{\!\!\!\!\!\!\!t}_{j,G,\gamma}}(k-k')\xi}
|\widehat{\psi}^{(G)}(2^{-\;D^{\!\!\!\!\!\!\!t}_{j,G,\gamma}}\;\xi)|^2\rho_{E_0}(\xi)^{-2H_0-d}\rmd\xi\;.
\]
Set $\zeta =2^{-\;D^{\!\!\!\!\!\!\!t}_{j,G,\gamma}}\;\xi$ and integrate by parts with respect to $\zeta_{\ell}$. Hence
\[
(k_{\ell}-k'_{\ell})\mathbb{E}(c_{j,G,\gamma}^{k}c_{j,G,\gamma}^{k'})=2^{j\mathrm{Tr}(D_{j,G,\gamma})}\int_{\mathbb{R}^{d}}(k_{\ell}-k'_{\ell})\rme^{\rmi(k-k')\zeta}
|\widehat{\psi}^{(G)}(\zeta)|^2\rho_{E_0}(2^{\;D^{\!\!\!\!\!\!\!t}_{j,G,\gamma}}\;\zeta)^{-2H_0-d}\rmd\zeta.
\]
Recall that the pseudo--norm may be assumed to be $\mathcal{C}^{\infty}(\mathbb{R}^{d}\setminus\{0\})$. Since the multi resolution analysis is $\mathcal{C}^1$
\[
\begin{array}{lll}
(k_{\ell}-k'_{\ell})\mathbb{E}(c_{j,G,\gamma}^{k}c_{j,G,\gamma}^{k'})&=&-2^{j\mathrm{Tr}(D_{j,G,\gamma})}\int_{\mathbb{R}^{d}}\rme^{\rmi(k-k')\zeta}
\frac{\partial}{\partial \zeta_{\ell}}\left(|\widehat{\psi}^{(G)}(\zeta)|^2\rho_{E_0}(2^{\;D^{\!\!\!\!\!\!\!t}_{j,G,\gamma}}\;\zeta)^{-2H_0-d}\right)\rmd\zeta\\
&=&-2^{j\mathrm{Tr}(D_{j,G,\gamma})}\int_{\mathbb{R}^{d}}\rme^{\rmi(k-k')\zeta}
\left(\frac{\partial}{\partial \zeta_{\ell}}|\widehat{\psi}^{(G)}(\zeta)|^2\right)\rho_{E_0}(2^{\;D^{\!\!\!\!\!\!\!t}_{j,G,\gamma}}\;\zeta)^{-2H_0-d}\rmd\zeta\\
&&-2^{jTr(D_{j,G,\gamma})}\int_{\mathbb{R}^{d}}\frac{\rme^{\rmi(k-k')\zeta}|\widehat{\psi}^{(G)}(\zeta)|^2}{\rho_{E_0}(2^{\;D^{\!\!\!\!\!\!\!t}_{j,G,\gamma}}\;\zeta)^{2H_0+d+1}}
\left(2^{\gamma_{\ell}}\frac{\partial}{\partial \zeta_{\ell}}(\rho_{E_0})(2^{\;D^{\!\!\!\!\!\!\!t}_{j,G,\gamma}}\;\zeta)\right)\rmd\zeta\;.
\end{array}
\]
An approach similar to the proof of inequalities~(\ref{eq:WC2}) yields
\begin{equation}\label{eq:ineq1}
2^{j\mathrm{Tr}(D_{j,G,\gamma})}\left|\int_{\mathbb{R}^{d}}e^{i(k-k')\zeta}
\left(\frac{\partial}{\partial \zeta_{\ell}}|\widehat{\psi}^{(G)}(\zeta)|^2\right)\rho_{E_0}(2^{\;D^{\!\!\!\!\!\!\!t}_{j,G,\gamma}}\;\zeta)^{-2H_0-d}d\zeta\right|\leq Cj^{d/\rho_{\min}(E_0)}2^{-2jH_0}\;.
\end{equation}
Further, differentiate the homogeneity relationship satisfied by $\rho_{\;E_0^{\!\!\!\!\!\!t}}\;$ and deduce that for any $a>0$ and $z\in\mathbb{R}^d$
\begin{equation}\label{eq:homgrad}
a^{\;E_0^{\!\!\!\!\!\!t}}\;(\overrightarrow{\mathrm{grad}}(\rho_{\;E_0^{\!\!\!\!\!\!t}}\;))(a^{\;E_0^{\!\!\!\!\!\!t}}\;z)=
a(\overrightarrow{\mathrm{grad}}(\rho_{\;E_0^{\!\!\!\!\!\!t}}\;))(z)\;.
\end{equation}
For any $y\in\mathbb{R}^{d}\setminus\{0\}$, let $r=|y|_{\;E_0^{\!\!\!\!\!\!t}}\;$. Then set $j=[\log_2(r)]$ and remark that $|\Theta|_{\;E_0^{\!\!\!\!\!\!t}}\;=|2^{-j\;E_0^{\!\!\!\!\!\!t}}\;y|_{\;E_0^{\!\!\!\!\!\!t}}\;\in [1/2,2]$ and hence that $\Theta$ belongs to the compact set $\mathcal{C}(1/2,2,\;E_0^{\!\!\!\!\!\!\!t})\;=\{\theta,\,|\theta|_{\;E_0^{\!\!\!\!\!\!t}}\;\in [1/2,2]\}$. Relationship~(\ref{eq:homgrad}) applied with $a=2^j$ and $z=\Theta$ then implies
\[
2^{\;D_{j,G,\gamma}^{\!\!\!\!\!\!t}}\;\overrightarrow{\mathrm{grad}}(\rho_{\;E_0^{\!\!\!\!\!\!t}}))(2^{j \;E_0^{\!\!\!\!\!\!t}}\;\Theta)=2^{-j \;E_0^{\!\!\!\!\!\!t}\;+\;D_{j,G,\gamma}^{\!\!\!\!\!\!t}}2^j(\overrightarrow{\mathrm{grad}(\rho_{\;E_0^{\!\!\!\!\!\!t}})}\;(\Theta))
\]
Take the norm of each member of the equality and deduce that for any $y\in\mathbb{R}^d\setminus\{0\}$ satisfying $j=[\log_2|y|_{\;E_0^{\!\!\!\!\!\!t}}\;]$
\[
|2^{\;D_{j,G,\gamma}^{\!\!\!\!\!\!t}}\overrightarrow{\mathrm{grad}}(\rho_{E_0^{\!\!\!\!\!\!t}}\;))(y)|\leq C 2^j |2^{-j \;E_0^{\!\!\!\!\!\!t}\;+\;D_{j,G,\gamma}^{\!\!\!\!\!\!t}}\;|
\]
where $C=\sup_{\Theta\in \mathcal{C}(1/2,2,\;E_0^{\!\!\!\!\!\!t})}\;|\overrightarrow{\mathrm{grad}}(\rho_{\;E_0^{\!\!\!\!\!\!t}})\;(\Theta)|$.\\
Lemma~2.1 of~\cite{5} and the definition of $j$ imply that
\[
|2^{\;D_{j,G,\gamma}^{\!\!\!\!\!\!t}}\;\overrightarrow{\mathrm{grad}}(\rho_{\;E_0^{\!\!\!\!\!\!t}}))(y)|\leq C 2^j |j|^{d/\rho_{\min}}\leq C |y|_{\;E_0^{\!\!\!\!\!\!t}}\;|\log(|y|_{\;E_0^{\!\!\!\!\!\!t}}\;)|^{d/\rho_{\min}(E_0)}\;.
\]
Set now $y=2^{\;D_{j,G,\gamma}^{\!\!\!\!\!\!t}}\;\zeta$. One has
\begin{eqnarray*}
|2^{\;D_{j,G,\gamma}^{\!\!\!\!\!\!\!t}}\;(\overrightarrow{\mathrm{grad}}(\rho_{\;E_0^{\!\!\!\!\!\!t}}\;))(2^{\;D_{j,G,\gamma}^{\!\!\!\!\!\!\!t}}\;\zeta)|
&\leq& C |2^{\;D_{j,G,\gamma}^{\!\!\!\!\!\!\!t}}\;\zeta|_{\;E_0^{\!\!\!\!\!\!t}}\;|\log(|2^{\;D_{j,G,\gamma}^{\!\!\!\!\!\!\!t}}\;\zeta|_{\;E_0^{\!\!\!\!\!\!t}}\;)|^{d/\rho_{\min}(E_0)}\\
&\leq& C 2^j (j+|\log(|\zeta|_{\;E_0^{\!\!\!\!\!\!t}}\;|)|)^{2d/\rho_{\min}(E_0)} |\zeta|_{\;E_0^{\!\!\!\!\!\!t}}\;.
\end{eqnarray*}
Since for any $\ell\in\{1,\cdots,d\}$
\[
2^{\gamma_{\ell}}\left|\left(\frac{\partial}{\partial \zeta_{\ell}}(\rho_{E_0^*})\right)(2^{\;D_{j,G,\gamma}^{\!\!\!\!\!\!t}}\;\zeta)\right|\leq \left|2^{\;D_{j,G,\gamma}^{\!\!\!\!\!\!t}}(\overrightarrow{\mathrm{grad}}(\rho_{E_0}))(2^{\;D_{j,G,\gamma}^{\!\!\!\!\!\!t}}\;\zeta)\right|
\]
it yields the following inequality
\begin{equation}\label{eq:ineq2}
\left|2^{j\mathrm{Tr}(D_{j,G,\gamma})}\int_{\mathbb{R}^{d}}\frac{\rme^{\rmi(k-k')\zeta}|\widehat{\psi}^{(G)}(\zeta)|^2}{\rho_{E_0}(2^{\;D_{j,G,\gamma}^{\!\!\!\!\!\!t}}\;\zeta)^{2H_0+d+1}}
\left(2^{\gamma_{\ell}}\frac{\partial}{\partial \zeta_{\ell}}(\rho_{E_0})(2^{\;D_{j,G,\gamma}^{\!\!\!\!\!\!t}}\;\zeta)\right)d\zeta\right|\leq 2^{-2j H_0} |j|^{2d/\rho_{\min}(E_0)}\;.
\end{equation}
Combining inequalities~(\ref{eq:ineq1}) and (\ref{eq:ineq2}) then yield inequality~(\ref{eq:WC1}).

Remark now that
\begin{lemma}\label{lem:estimnj}
Let $D_0$ an admissible diagonal anisotropy satisfying
$\mathrm{Tr}(D_0)=d$. Recall that $\Gamma_j(D_0)$ is defined
by~(\ref{e:defGammajD}). There exists some $C_1,C_2>0$ such that
\[
C_1 j^d 2^{jd}\leq \mathrm{card}(\Gamma_j(D_0))\leq j^d 2^{jd}\;.
\]
\end{lemma}
{\bf Proof.} Indeed, since the norms $|\cdot|_{\ell_1}$ and $|\cdot|_{\ell^{\infty}}$ on $\mathbb{R}^d$ are equivalent, there exists some $C_1,C_2>0$ such that
\[
C_1\max_{\ell}|k_{\ell}|^{1/\lambda_{\ell}}\leq |k|_{D_0}\leq C_2\max_{\ell}|k_{\ell}|^{1/\lambda_{\ell}}\;.
\]
The conclusion follows since it is quite clear since that
\[
\mathrm{card}\{k,\max_{\ell}|k_{\ell}|^{1/\lambda_{\ell}}\leq j2^j\}=\mathrm{card}\{k,\max_{\ell}|k_{\ell}|\leq j^{\lambda_{\ell}}2^{j\lambda_{\ell}}\}=\prod_{\ell}\left(j^{\lambda_{\ell}}2^{j\lambda_{\ell}}\right)=j^d2^{d}\;,
\]
using the fact that $\lambda_1+\cdots+\lambda_{\ell}=d$.\\

The proof of Proposition~\ref{pro:X1} is then based on the two following results which are a slight modification of Theorem II.1 and II.7 of~\cite{10}. We recall the proofs for completeness.\\
We denote
\[
c_p=\mathbb{E}(|g^0_{j,G,\gamma}|^p)\;.
\]
where $(g_{j,G,\gamma}^{k})$ is the stationary Gaussian sequence of the normalized
wavelet coefficients defined by~(\ref{eq:defgjk}). Since, by
Lemma~\ref{LemKey}, the wavelet coefficients are weakly dependent,
we can state a central limit theorem for the sequence
$(g_{j,G,\gamma}^{k})_{j\in\mathbb{N},(G,\gamma)\in
I_j,k\in\Gamma_j}$ which is a slight modified version of
Lemma~II.4 of \cite{10}
\begin{lemma}\label{lem:CRK1}
Let $p\in (1,+\infty)$ and $(g_{j,G,\gamma}^{k})$ the Gaussian
sequence defined by~(\ref{eq:defgjk}). Set
$n_j=\mathrm{Card}(\Gamma_j(D_0))$. Then almost surely when
$j\to\infty$
\[
n_j^{-1}\left(\sum_{k\in
\Gamma_j(D_0)}|g_{j,G,\gamma}|^p\right)\to c_p\;.
\]
\end{lemma}
{\bf Proof.} By Lemma~\ref{LemKey} the sequence $(g_{j,G,\gamma}^{k})$ is weakly correlated in the sense of~\cite{10}--that is satisfies the assumption (H) of \cite{10}. We follow the main line of~\cite{10} and first give an upper bound of
\[
\mathbb{E}\left|\sum_{k\in \Gamma_j}\left(|g_{j,G,\gamma}|^p-c_p\right)\right|^2\;.
\]
Using the same approach that in~\cite{10} (see Lemma II.3) we get that
\[
\mathbb{E}\left|\sum_{k\in \Gamma_j}\left(|g_{j,G,\gamma}|^p-c_p\right)\right|^2\leq C_j c_{2p}\sum_{(k,k')\in \Gamma_j^2}\frac{1}{(1+|k-k'|)^2}\;,
\]
with $C_j=j^{2d/\rho_{\min}(E_0)}$ by weak correlation of the wavelet coefficients. Set $\ell=k-k'$. Hence
\[
\sum_{(k,k')\in \Gamma_j^2}\frac{1}{(1+|k-k'|)^2}\leq \sum_{k\in \Gamma_j}\sum_{\ell\in 2.\Gamma_j}\frac{1}{(1+|\ell|)^2}\leq C_j j2^{j}\sum_{\ell\in 2.\Gamma_j}\frac{1}{(1+|\ell|)^2}\;.
\]
Remark now that
\[
\sum_{\ell\in 2.\Gamma_j}\frac{1}{(1+|\ell|)^2}\leq \sum_{\ell\in 2.\Gamma_j}\frac{1}{(1+|\ell|_{D_0})^{2/\rho_{\max}}}\leq j^{d-\delta}2^{j(d-\delta)}\;,
\]
with $\delta=2/\rho_{\max}(E_0)>0$ by comparison with an integral and Proposition~2.3 in~\cite{6}.\\
Thereafter the end of the proof is exactly the same that in Theorem II.1 in~\cite{10}.\\

In an analogous way, one can give a result on the asymptotic behavior of
\[
\frac{1}{\sqrt{|\log(n_j)|}}\left(\max_{k\in \Gamma_j}|g_{j,G,\gamma}|\right)\;.
\]
\begin{lemma}\label{lem:CRK2}
Almost surely
\[
0<\liminf_{j\to\infty}\frac{1}{\sqrt{|\log(n_j)|}}\left(\max_{k\in \Gamma_j}|g_{j,G,\gamma}|\right)\leq \limsup_{j\to\infty}\frac{1}{\sqrt{|\log(n_j)|}}\left(\max_{k\in \Gamma_j}|g_{j,G,\gamma}|\right)<\infty\;.
\]
\end{lemma}
{\bf Proof.} The proof is exactly the same than these of Lemmas~II.8 and~II.10 in~\cite{10}.



 { \footnotesize
\medskip
\medskip
 \vspace*{1mm}

\noindent {\it M. Clausel}\\
Laboratoire Jean Kuntzmann, Universit\'e de Grenoble, CNRS, F38041 Grenoble Cedex 9\\
E-mail: {\tt marianne.clausel@imag.fr}\\

\noindent {\it B.Vedel}\\
Laboratoire de Mathematiques et Applications des Math\'ematiques,\\ Universite de Bretagne Sud, Universit\'e Europ\'eene de Bretagne
Centre Yves Coppens, Bat. B, 1er et., Campus de Tohannic BP 573,
56017 Vannes, France.\\
E-mail: {\tt vedel@univ-ubs.fr }

 \end{document}